\documentclass[11pt]{article}

\usepackage{IEEEtrantools}
\usepackage{latexsym}
\usepackage{amsfonts}
\usepackage{amssymb}
\usepackage{amsmath,amsfonts,amsthm}
\usepackage{mathrsfs}
\usepackage{enumerate}
\usepackage{dsfont}

\usepackage{tikz-cd}

\usepackage{multirow}
\usepackage{tikz}
\usetikzlibrary{arrows,decorations.markings}
\usepackage{bookmark}
\usepackage{hyperref}
\hypersetup{
     colorlinks   = true,
     citecolor    = blue
}

\newcommand{\bbbt}{\mathbb{T}}
\newcommand{\scrt}{\mathscr{T}}

\newcommand{\be}{\begin{equation}}
\newcommand{\ee}{\end{equation}}
\newcommand{\bea}{\begin{eqnarray}}
\newcommand{\eea}{\end{eqnarray}}
\newcommand{\bean}{\begin{eqnarray*}}
\newcommand{\eean}{\end{eqnarray*}}
\newcommand{\brray}{\begin{array}}
\newcommand{\erray}{\end{array}}
\newcommand{\biearray}{\begin{IEEEarray}{rCl}}
\newcommand{\eiearray}{\end{IEEEarray}}

\newcommand{\newsection}[1]{\setcounter{equation}{0}
\setcounter{dfn}{0}
\section{#1}}

\newtheorem{dfn}{Definition}[section]
\newtheorem{thm}[dfn]{Theorem}
\newtheorem{lmma}[dfn]{Lemma}
\newtheorem{ppsn}[dfn]{Proposition}
\newtheorem{crlre}[dfn]{Corollary}
\newtheorem{xmpl}[dfn]{Example}
\newtheorem{rmrk}[dfn]{Remark}
\newcommand{\restr}[1]{|_{#1}}

\newcommand{\bdfn}{\begin{dfn}\rm}
\newcommand{\bthm}{\begin{thm}}
\newcommand{\blmma}{\begin{lmma}}
\newcommand{\bppsn}{\begin{ppsn}}
\newcommand{\bcrlre}{\begin{crlre}}
\newcommand{\bxmpl}{\begin{xmpl}}
\newcommand{\brmrk}{\begin{rmrk}\rm}

\newcommand{\edfn}{\end{dfn}}
\newcommand{\ethm}{\end{thm}}
\newcommand{\elmma}{\end{lmma}}
\newcommand{\eppsn}{\end{ppsn}}
\newcommand{\ecrlre}{\end{crlre}}
\newcommand{\exmpl}{\end{xmpl}}
\newcommand{\ermrk}{\end{rmrk}}

\newcommand{\bbc}{\mathbb{C}}
\newcommand{\bbz}{\mathbb{Z}}
\newcommand{\bbn}{\mathbb{N}}

\newcommand{\cla}{\mathcal{A}}
\newcommand{\clb}{\mathcal{B}}
\newcommand{\clh}{\mathcal{H}}

\newcommand{\clk}{\mathcal{K}}

\newcommand{\prf}{\noindent{\it Proof\/}: }

\def \qed { \mbox{}\hfill
$\Box$\vspace{1ex}}

\begin{document}

%%%%%%%%%%%%%%%%%%%%%%%%%%%%%%%%%
%%%%%%%%%%%%%%%%%%%%%%%%%%%%%%%%%

\author{\sc{Satyajit Guin, Bipul Saurabh}}
\title{Equivariant spectral triple for the quantum group $U_q(2)$ for complex deformation parameters}
\maketitle

%%%%%%%%%%%%%%%%%%%%%%%%%%%%%%%%%%
%%%%%%%%%%%   ABSTRACT    %%%%%%%%%%%%%%%%
%%%%%%%%%%%%%%%%%%%%%%%%%%%%%%%%%%

\begin{abstract}
Let $q=|q|e^{i\pi\theta}$ be a nonzero complex number such that $|q|\neq 1$ and consider the compact quantum group $U_q(2)$. For $\theta\notin\mathbb{Q}\setminus\{0,1\}$, we obtain the $K$-theory of the underlying $C^*$-algebra $C(U_q(2))$. We construct a spectral triple on $U_q(2)$ which is equivariant under its own comultiplication action. The spectral triple obtained 
 here is even, $4^+$-summable, non-degenerate, and the Dirac operator acts on two copies of the $L^2$-space of $U_q(2)$. The $K$-homology class of the associated Fredholm module is shown to be nontrivial.
\end{abstract}
\bigskip

{\bf AMS Subject Classification No.:} {\large 58}B{\large 34}, {\large 58}B{\large 32}, {\large46}L{\large 87}.

{\bf Keywords.} Compact quantum group, spectral triple, quantum unitary group, equivariance.
\bigskip

%%%%%%%%%%%%%%%%%%%%%%%%%%%%%%%%%%
%%%%%%%%%   INTRODUCTION    %%%%%%%%%%%%%%%
%%%%%%%%%%%%%%%%%%%%%%%%%%%%%%%%%%

\newsection{Introduction}

In \cite{KasMeyRoyWor-2016aa}, Woronowicz et al. defined a family of $q$-deformations of $SU(2)$ for $q\in\mathbb{C}\setminus\{0\}$. This agrees with the compact quantum group $SU_q(2)$ (\cite{Wor1},\cite{Wor2}) when $q$ is real. But for $q\in\mathbb{C}\setminus\mathbb{R},\,SU_q(2)$ is not a compact quantum group, rather a braided compact quantum group in a suitable tensor category. In \cite{MeyRoyWor-2016aa}, Woronowicz et al. have shown that for a compact quantum group $\mathbb{G}=(\mathbb{A},\Delta)$ and a braided compact quantum group $\mathbb{B}$ over $\mathbb{G}$, the semidirect product $\mathbb{A}\boxtimes\mathbb{B}$ becomes a compact quantum group. This semidirect product construction is a $C^*$-algebraic analogue of ``bosonisation'' defined by Majid \cite{Maj-1994aa}. Taking $\mathbb{B}=SU_q(2)$ for $q\in\mathbb{C}\setminus\{0\}$ and $\mathbb{A}=C(\mathbb{T})$, we obtain a genuine compact quantum group. In (Sec. $6$, \cite{KasMeyRoyWor-2016aa}), it is shown that this compact quantum group is the coopposite of the compact quantum group $U_q(2)$ defined in \cite{ZhaZha-2005aa}. In \cite{GuiSau-2020aa}, authors have analysed the quantum group structure of $U_q(2)$ by describing all its finite dimensional irreducible unitary representations up to equivalence. The matrix coefficients are expressed in terms of the generators of $C(U_q(2))$ using the little $q$-Jacobi polynomials, their norms have been computed, and the explicit form of the Peter-Weyl decomposition is obtained. The representation ring structure is also determined by showing how the tensor product of two irreducible representations decomposes into irreducible components. These works lay the foundation for investigating geometrical aspects of $U_q(2)$ in the sense of Connes (\cite{Con-1994aa},\cite{Con-1996aa},\cite{CM-1995aa}), and purpose of the present article is precisely that.

Central notion in Connes' formulation of Noncommutative geometry is the notion of ``Spectral triple''. Some nice properties e.g. finite-summability, nontriviality on the spectral triple are also desirable. Moreover, if a classical or quantum group acts on the underlying space then we desire to get a spectral triple that is equivariant under the group action. A quantum group $G$ acts on itself via its comultiplaction action. Two natural questions arise here namely, whether there exists a finitely-summable nontrivial 
spectral triple for $G$ that is equivariant under the group action, and if it exists, then can we construct one such concretely? 
In (\cite{CP-2003aa},\cite{DLSSV-2005aa}), both the problems for the quantum group $SU_q(2)$ are settled. Neshveyev-Tuset (\cite{NesTus-2010aa},\cite{NesTus-2010ab}) answered the first problem for the $q$-deformation $G_q$ of a simply connected semisimple compact Lie group $G$. The Dirac operator obtained for $G_q$ is isospectral to its classical Dirac operator of $G$. However, this construction is existential in nature as it involves twisting the classical Dirac operator by a Drinfeld unitary whose concrete realisation on a Hilbert space is difficult to obtain. As a consequence, showing regularity or exhibiting the local index formula \cite{CM-1995aa} for these geometric spaces seems to be a difficult task. Our approach is mainly along the line of (\cite{CP-2003aa},\cite{DLSSV-2005aa}) which primarily relies upon analysing the matrix elements of the irreducible representation of $U_q(2)$. Since the Dirac operators for $SU_q(2)$ in this approach is concretely realised, the local index formula is explicitly obtained (\cite{Con-2004aa},\cite{DLSSV-2006aa}). This approach also works for homogeneous spaces of $U_q(2)$ which we have discussed in \cite{GS}. 

We now briefly discuss our work here. For $q=|q|e^{\sqrt{-1}\pi\theta}$, we assume that $\theta\notin\mathbb{Q}\setminus\{0,1\}$. Note that for $\theta=0,1$, the deformation parameter $q$ is real. As in \cite{ZhaZha-2005aa}, we denote by $a,b,D$ the generators of the underlying $C^*$-algebra $C(U_q(2))$ and by $h$ its Haar state. It is known that $h$ is faithful (Thm. $2.8$ in \cite{GuiSau-2020aa}). The GNS representation of $C(U_q(2))$ on $L^2(h)$ is denoted by $\pi_h$. In Sec. \ref{Sec 2}, we briefly recall the notion of equivariant nontrivial spectral triple and the compact quantum group $U_q(2)$. In Sec. \ref{Sec 3}, we compute the $K$-theory of $C(U_q(2))$ along with the generators. This is achieved by showing that certain closed ideal in $C(U_q(2))$ can be identified with $\mathcal{K}(\ell^2(\bbn))\otimes\mathbb{A}_\theta$, where $\mathbb{A}_\theta$ is the noncommutative torus and $\mathcal{K}(\ell^2(\bbn))$ denotes the space of compact operators acting on $\ell^2(\bbn)$. It turns out that one of the generators, say $P_\theta$, of $K_0(C(U_q(2)))$ is same as the Powers-Rieffel projection in $\mathbb{A}_\theta$ if one replaces $Pb$ and $PD$, where $P$ denotes the element $\mathds{1}_{\{bb^*=1\}}$ in $C(U_q(2))$, with the standard unitary generators of $\mathbb{A}_\theta$. In Sec. \ref{Sec 4}, we compute the action of the generators of $U_q(2)$ on the orthonormal basis elements in $L^2(h)$.  Sec. \ref{Sec 5} deals with the fixed point subspace, say $E_1$, under the action of $bb^*$. Observe that $\pi_h(P)$ is the orthogonal projection onto $E_1$. Sec. \ref{Sec 6} is the heart of the article where we construct a $4^+$-summable even spectral triple on $L^2(h)\otimes\mathbb{C}^2$ which is equivariant under the comultiplaction action of $U_q(2)$. The Dirac operator is shown to be non-degenerate, which means that our Dirac operator is really a Dirac operator for the full tangent bundle rather than that of some lower dimensional subbundle. In Sec. \ref{Sec 7}, the $K$-homology class of the associated Fredholm module is shown to be nontrivial. This is achieved by pairing it with the class of projection  $P_\theta$ and the required Fredholm index turns out to be nonzero. Finally in Sec. \ref{Sec 8}, we have computed the spectral dimension \cite{CP-2016aa} of $U_q(2)$.

Throughout the paper the letter `$i$' will always denote an index, and whenever the complex indeterminate appears we denote it by $\sqrt{-1}\,$.

%%%%%%%%%%%%%%%%%%%%%%%%%%%%%%%%%%%%%%%%%
%%%%%%%%%%%%%%%     Preliminaries    %%%%%%%%%%%%%%%%%%
%%%%%%%%%%%%%%%%%%%%%%%%%%%%%%%%%%%%%%%%%

\newsection{Preliminaries}\label{Sec 2}

This section exhibits some basic notions and notations in use. Although well-known but we briefly discuss nontrivial spectral triple and equivariant Dirac operator here for the reader's convenience. The compact quantum group $U_q(2)$ is also briefly described and for detail discussion on it we refer to \cite{GuiSau-2020aa}. 
\medskip

\noindent\textbf{Nontrivial spectral Triple:} 
Let $\mathbb{A}$ be a unital $C^*$-algebra and $\cla\subseteq\mathbb{A}$ be a unital dense $\star$-subalgebra. An {\it odd} spectral triple is a tuple $(\cla,\clh,\mathcal{D})$, where $\clh$ is a separable Hilbert space on which $\cla$ acts as bounded operators, $\mathcal{D}$ is an unbounded self-adjoint operator with compact resolvent such that the commutator $[\mathcal{D},a]$ extends to a bounded operator on $\clh$ for all $a\in\cla$. If there is a $\bbz_2$-grading operator $\gamma$ acting on $\clh$ such that it commutes with $a$ for all $a\in\cla$ and anticommutes with $\mathcal{D}$, then $(\cla,\clh,\mathcal{D},\gamma)$ is called an {\it even} spectral triple. The operator $\mathcal{D}$ is widely referred as the Dirac operator. Since the kernel of $\mathcal{D}$ is finite dimensional, without loss of generality, one can assume that $\mathcal{D}$ has trivial kernel (P. $316$ in \cite{Con-1994aa}, P. $446$ in \cite{GVF-2001aa}). If for $p>0,\,|\mathcal{D}|^{-p}$ lies in the Dixmier ideal $\mathcal{L}^{(1,\infty)}\subseteq\mathcal{B}(\clh)$, then we say that the spectral triple is $p^+$-summable \cite{Con-1994aa}. Note that given $\mathcal{A}\subseteq\mathbb{A}$, there is no general recipe to construct a Dirac operator on $\mathcal{A}$ satisfying certain desirable properties. If $\cla$ is countably generated in $\mathbb{A}$, then existence of spectral triple $(\cla,\clh,\mathcal{D})$ is known (\cite{BJ-1983aa}, Chap. $4$ in \cite{Con-1994aa}). However, if we demand further properties like finite summability, then given a dense subalgebra of a $C^*$-algebra it may not admit any finitely summable spectral triple \cite{Con-1989aa}.

An even spectral triple $(\cla,\clh,\mathcal{D},\gamma)$ induces a $\mathcal{K}$-homology class $[(\cla,\clh,F,\gamma)]$, where $F=\mathcal{D}|\mathcal{D}|^{-1}$, in $K^0(\mathbb{A})$ consisting of even Fredholm modules. To check nontriviality of this class, one pairs it with $K_0(\mathbb{A})$ through the Kasparov product, which we describe briefly. Given a projection $P\in M_n(\mathbb{A})$, define $\clh_n=\clh\otimes \bbc^n,\,\gamma_n=\gamma\otimes I_n,\,F_n=F\otimes I_n,\,P^+=\frac{1+\gamma_n}{2}P$, and $P^-=\frac{1-\gamma_n}{2}P$. Then, $\clh_n$ decomposes as $\clh_n^+\oplus\clh_n^-$ under the grading operator $\gamma_n$ and the operator $P^-F_nP^+:P^+\clh_n^+\longrightarrow P^-\clh_n^-$ is a Fredholm operator. Index of this Fredholm operator is the value of the $K_0-K^0$ pairing $\langle[P],[(\cla,\clh,F,\gamma)]\rangle$. An even spectral triple $(\clh,\mathcal{D},\gamma)$ on $\cla$ is called nontrivial if $\langle[P],[(\cla,\clh,F,\gamma)]\rangle$ is nonzero for some $[P]\in K_0(\mathbb{A})$. The requirement of nontrivial pairing is very crucial and reason is discussed in detail in \cite{CP-2003aa}.
\medskip

\noindent\textbf{Equivariant Dirac Operator:} We first recall the definition in the general setting from \cite{CP-2008aa}, and then discuss the situation relevant to this article.

Assume that a compact quantum group $\mathbb{G}$ has an action on $\mathbb{A}$ given by $\tau:\mathbb{A}\longrightarrow\mathbb{A}\otimes C(\mathbb{G})$ so that $(\mbox{id}\otimes\Delta)\tau=(\tau\otimes\mbox{id})\tau,\,\Delta$ being the coproduct. We have a $C^*$-dynamical system $(\mathbb{A},\mathbb{G},\tau)$. A covariant representation $(\pi,\mathbb{U})$ of the $C^*$-dynamical system $(\mathbb{A},\mathbb{G},\tau)$ consists of a unital $\star$-representation $\pi:\mathbb{A}\longrightarrow\mathbb{B}(\mathcal{H})$ and a unitary representation $\mathbb{U}$ of $\mathbb{G}$ on $\mathcal{H}$, i,e. a unitary element of the multiplier algebra $M(\mathcal{K}(\mathcal{H})\otimes C(\mathbb{G}))$ such that $(\pi\otimes\mbox{id})\tau(a)=\mathbb{U}(\pi(a)\otimes I)\mathbb{U}^*$ for all $a\in\mathbb{A}$. A $\mathbb{G}$-equivariant odd spectral triple for $(\mathbb{A},\mathbb{G},\tau)$ is a quadruple $(\mathcal{H},\pi,\mathbb{U},\mathcal{D})$ such that the following conditions are satisfied~:
\begin{enumerate}[(i)] 
\item $(\pi,\mathbb{U})$ is a covariant representation of $(\mathbb{A},\mathbb{G},\tau)$ on $\mathcal{H}$,
\item $\pi$ is faithful,
\item $\mathbb{U}(\mathcal{D}\otimes I)\mathbb{U}^*=\mathcal{D}\otimes I$,
\item $(\mathcal{H},\pi,\mathcal{D})$ is an odd spectral triple for $\mathbb{A}$.
\end{enumerate} 
Similarly, a $\mathbb{G}$-equivariant even spectral triple for $(\mathbb{A},\mathbb{G},\tau)$ is a tuple $(\mathcal{H},\pi,\mathbb{U},\mathcal{D},\gamma)$, where $\gamma$ is a grading operator such that $(\mathcal{H},\pi,\mathbb{U},\mathcal{D})$ satisfies all the above conditions and $\mathbb{U}(\gamma \otimes I)\mathbb{U}^*=\gamma\otimes I$. If all the notations are clear from the context, then we simply say that the Diarc operator $\mathcal{D}$ is $\mathbb{G}$-equivariant. 

Now, we discuss the following special case which is relevant to our situation in this article. Suppose $\mathbb{G}=(C(\mathbb{G}),\Delta)$ is a compact quantum group with a faithful Haar state $h$. Then, $\mathbb{G}$ acts on itself via the comultiplication action $\Delta$. The invariant state of this action will be the Haar state $h$. Let $(u^{\alpha})_{\alpha \in \widehat{\mathbb{G}}}$ be a complete family of irreducible unitary corepresentations of $\mathbb{G}$. Suppose $u^{\alpha}$ acts on the Hilbert space $H_{\alpha}$ with dimension $N_{\alpha}$. Using an orthonormal basis of $H_{\alpha}$, one can identify $u^{\alpha}$ with a  matrix $(u_{ij}^{\alpha})_{1\leq i,j\leq N_{\alpha}}$, where $u_{ij}^{\alpha}$ are elements of $C(\mathbb{G})$ such that the following holds~:
\begin{IEEEeqnarray}{rCl}
\sum_m u_{im}^{\alpha}(u_{km}^{\alpha})^*&=&\delta_{ik}\,, \\
\sum_m u_{mj}^{\alpha}(u_{mk}^{\alpha})^*&=&\delta_{jk}\,, \\
\Delta(u_{ij}^{\alpha}) &=& \sum_k u_{ik}^{\alpha}\otimes u_{kj}^{\alpha}\,.
\end{IEEEeqnarray}
Moreover, one has  the following results (see \cite{Wor1}),
\begin{enumerate}[(i)]
\item Let $\mathcal{O}(\mathbb{G})$ be the linear span of $(u^{\alpha})_{\alpha \in \widehat{\mathbb{G}}}$. Then, this is a dense $\star$-subalgebra of $C(\mathbb{G})$.
\item The set  $\left\{e_{ij}^{\alpha}:=\frac{u_{ij}^{\alpha}}{\|u_{ij}^{\alpha}\|}:1\leq i,j\leq N_{\alpha}, \alpha \in \widehat{\mathbb{G}}\right\}$ 
is an orthonormal basis of $L^2(h)$.
\item For $\alpha \in \widehat{\mathbb{G}}$, one has $\|u_{ij}^{\alpha}\|=\|u_{i^{'}j^{'}}^{\alpha}\|$ for any $1 \leq i,i^{'},j,j^{'}\leq N_{\alpha}$.
\end{enumerate}
Let $\pi:C(\mathbb{G})\longrightarrow \clb(L^2(h))$ be the GNS representation associated to the faithful Haar state $h$, and $\mathbb{U}:L^2(h)\longrightarrow L^2(h)\otimes C(\mathbb{G})$ be the unitary corepresentation of $\mathbb{G}$ implemented by the comultiplication $\Delta$. Then, it is not difficult to check that $(\pi,\mathbb{U})$ is a covariant representation of the $C^*$-dynamical system $(C(\mathbb{G}),\mathbb{G},\Delta)$. Given any bounded linear functional $\rho$ on $C(\mathbb{G})$, define $\mathbb{U}_{\rho}=(\mbox{id}\otimes\rho)\mathbb{U}$. Observe that $\mathbb{U}_{\rho}\in\mathcal{B}(L^2(h))$ and $\mathbb{U}_{\rho}$ takes $\mathcal{O}(\mathbb{G})$ to $\mathcal{O}(\mathbb{G})$. In this set-up, a Dirac operator $\mathcal{D}$ is equivariant if it has domain $\mathcal{O}(\mathbb{G})$ and the following condition is satisfied~:
\begin{IEEEeqnarray}{rCl}
\mathbb{U}_{\rho}\mathcal{D}(x)=\mathcal{D}\mathbb{U}_{\rho}(x)\quad\mbox{ for all }\,x\in\mathcal{O}(\mathbb{G})
\end{IEEEeqnarray}
and for any bounded linear functional $\rho$ on $C(\mathbb{G})$.

\bppsn\label{equivariance}
Let $(C(\mathbb{G}), \pi, T)$ be an odd spectral triple such that $Tu_{ij}^{\alpha}=d(\alpha,i)u_{ij}^{\alpha}$ for $1\leq i,j\leq N_{\alpha},\, \alpha \in \widehat{\mathbb{G}}$. Then, one has the followings~:
\begin{enumerate}[$(i)$]
\item The spectral triple $(C(\mathbb{G}), \pi, T)$ is equivariant under the comultiplication action $\Delta$ of $\mathbb{G}$.
\item The tuple $(\pi\oplus \pi, \mathbb{U}\oplus \mathbb{U})$ is a covariant representation of the $C^*$-dynamical system $(C(\mathbb{G}),\mathbb{G},\Delta)$ on the Hilbert space $L^2(h)\oplus L^2(h)$. Moreover, the following tuple
\[
\left(C(\mathbb{G})\,,\, \pi\oplus \pi\,,\,  \mathscr{D}:=\left( {\begin{matrix}
   0  & T^*\\
  T & 0 \\ 
  \end{matrix} } \right)\,,\, \gamma:= \left( {\begin{matrix}
   1 & 0\\
  0 & -1 \\
  \end{matrix} } \right)\right)\]
is an even spectral triple that is equivariant under the comultiplication action $\Delta$ of $\mathbb{G}$. 
\end{enumerate}
\eppsn
\prf To prove part $(i)$, take a bounded linear functional $\rho$ of $C(\mathbb{G})$. 
 It suffices to show  that $T\mathbb{U}_{\rho}(u_{ij}^{\alpha})=\mathbb{U}_{\rho}T(u_{ij}^{\alpha})$ for all $1\leq i,j\leq N_{\alpha},\, \alpha \in \widehat{\mathbb{G}}$. Observe the following, 
\begin{IEEEeqnarray*}{rCl}
 T\mathbb{U}_{\rho}(u_{ij}^{\alpha})&=& \sum_k T(u_{ik}^{\alpha})\rho(u_{kj}^{\alpha}), \\
 &=&\sum_k d(\alpha, i) \rho(u_{kj}^{\alpha})u_{ik}^{\alpha} ,\\
 &=&d(\alpha, i)\mathbb{U}_{\rho}(u_{ij}^{\alpha}),\\
 &=&\mathbb{U}_{\rho}(T(u_{ij}^{\alpha}))\,.
\end{IEEEeqnarray*}
The first claim of part $(ii)$ follows from a straightforward verification. The remaining part follows from the following,
\begin{IEEEeqnarray*}{rCl}
 \left[\left( {\begin{matrix}
   0  & T^*\\
  T & 0 \\ 
  \end{matrix} } \right), \left( {\begin{matrix}
   \pi(a)  & 0\\
  0 & \pi(a) \\ 
  \end{matrix} } \right)\right]&=&\left( {\begin{matrix}
   0  & [T^*, \pi(a)]\\
  [T,\pi(a)] & 0 \\ 
  \end{matrix} } \right)
  = \left( {\begin{matrix}
   0  & -[T, \pi(a^*)]^*\\
  [T,\pi(a)] & 0 \\ 
  \end{matrix} } \right)\\
 \end{IEEEeqnarray*} 
 and
  \begin{IEEEeqnarray*}{rCl}
  \left( {\begin{matrix}
   \mathbb{U}  & 0\\
  0 &  \mathbb{U} \\ 
  \end{matrix} } \right) \left( {\begin{matrix}
   0  & T^*\otimes 1\\
  T \otimes 1 & 0 \\ 
  \end{matrix} } \right)\left( {\begin{matrix}
    \mathbb{U}^*   & 0\\
  0 &  \mathbb{U}^*  \\ 
  \end{matrix} } \right)&=&\left( {\begin{matrix}
   0  &  \mathbb{U} (T^*\otimes 1) \mathbb{U}^*\\
   \mathbb{U} (T \otimes 1)  \mathbb{U}^* & 0 \\ 
  \end{matrix} } \right)
  = \left( {\begin{matrix}
   0  &   T^*\otimes 1 \\
  T\otimes 1  & 0 \\ 
  \end{matrix} } \right)
\end{IEEEeqnarray*}
as $\mathbb{U}(T\otimes 1)\mathbb{U}^*=T\otimes 1$ from part $(i)$. This completes the proof.\qed
\medskip

\noindent\textbf{Compact quantum group $U_q(2)$:} 
Let us recall the definition of $U_q(2)$ from \cite{ZhaZha-2005aa}. For any nonzero complex number $q$, 
the $C^*$-algebra $C(U_q(2))$ is the universal $C^*$-algebra generated by $a,b,D$ satisfying the following relations~:
\begin{IEEEeqnarray}{lCl} \label{relations}
ba&=& q ab, \qquad a^*b=qba^*, \qquad   \qquad \qquad bb^*=b^*b, \qquad \qquad aa^*+bb^*=1, \nonumber \\
aD&=&Da, \qquad bD=q^2|q|^{-2}Db,  \qquad DD^*=D^*D=1,  \qquad a^*a+|q|^2b^*b=1.  
\end{IEEEeqnarray}
The compact quantum group structure is given by the following comultiplication,
\begin{IEEEeqnarray}{lCl} \label{comul}
\Delta(a)=a \otimes a-\bar{q}b \otimes Db^*\quad,\quad\Delta(b)=a\otimes b+b\otimes Da^*\quad,\quad\Delta(D)=D \otimes D.
\end{IEEEeqnarray}
Let $\mathscr{A}_q:=\mathcal{O}(U_q(2))$ be the $\star$-subalgebra of the $C^*$-algebra $C(U_q(2))$ generated by $a,b$ and $D$. The Hopf $\star$-algebra structure on it is given by the following~: 
\begin{IEEEeqnarray*}{lCl}
\mbox{antipode:} \quad S(a)=a^*,\,\,S(b)=-qbD^*,\,\,S(D)=D^*,\,\,S(a^*)=a,\,\,S(b^*)=-(\bar{q})^{-1}b^*D\,,\\ 
\mbox{counit:} \qquad\epsilon(a)=1,\,\,\epsilon(b)=0,\,\,\epsilon(D)=1\,.
\end{IEEEeqnarray*}
For $|q|\neq 1$, it is a compact quantum group of {\it non-Kac} type, whereas for $|q|=1$ it is of {\it Kac} type. In this article, we restrict ourselves to the case of $\,|q|\neq 1$. Faithful $C^*$-representations for the cases of $|q|\neq 1$ and $|q|=1$ are different and lie on different Hilbert spaces (Propn. $2.1$ and Propn. $6.1$ in \cite{GuiSau-2020aa}), and hence the case of $|q|=1$ may require different treatment. We plan to discuss this case in detail elsewhere. The classification obtained in (Thm. $7.1$, \cite{GuiSau-2020aa}) justifies that for the case of $|q|\neq 1$, it is enough to assume that $|q|<1$ in the context of this article, and henceforth we do so.

%%%%%%%%%%%%%%%%%%%%%%%%%%%%%%%%%%%%%%%%%%%%
%%%%%%%%%%%%%%%   K-groups   %%%%%%%%%%%%%%%%%%%%%%%
%%%%%%%%%%%%%%%%%%%%%%%%%%%%%%%%%%%%%%%%%%%%

\newsection{The \texorpdfstring{$K$}{}-theory of \texorpdfstring{$C(U_q(2))$}{}}\label{Sec 3}

Let $q$ be a nonzero complex number with the polar decomposition $q=|q|e^{\sqrt{-1}\pi\theta},\,\theta\in(-1,1]$, and $|q|\neq 1$. Note that for $\theta\neq 0,1\,$ the deformation parameter $q$ is non-real. The classification obtained in (Thm. $7.1$, \cite{GuiSau-2020aa}) justifies that it is enough to assume $|q|<1$ in the context of this article, and henceforth we do so. We further assume that $\theta\notin\mathbb{Q}\setminus\{0,1\}$. Let $\clh$ be the Hilbert space $\ell^2(\bbn)\otimes\ell^2(\bbz)\otimes\ell^2(\bbz)$. Consider the right shift operator $V:e_n\longmapsto e_{n+1}$ acting on $\ell^2(\bbn)$ and the unitary shift operator $U:e_k\longmapsto e_{k+1}$ acting on $\ell^2(\bbz)$. Let $N:e_n\longmapsto ne_n$ be the number operator. It is known that for $|q|<1$ the following representation of $C(U_q(2))$ on $\clh$~:
\begin{IEEEeqnarray}{lCl}\label{representation}
\pi(a)= \sqrt{1-|q|^{2N}}\,V \otimes 1 \otimes 1, \quad 
\pi(b) =q^N \otimes U \otimes 1, \quad
\pi(D) = 1 \otimes e^{-2\sqrt{-1}\pi\theta N}\otimes U
\end{IEEEeqnarray}
is faithful (Propn. $2.1$, \cite{GuiSau-2020aa}). We often omit the representation symbol $\pi$ for notational brevity.

We first deal with the case when $\theta$ is irrational, and mention at the end for $\theta=0,1$. For the standard orthonormal basis $\{e_n:n=0,1,\ldots\}$ of $\ell^2(\mathbb{N})$, we use the bra-ket notation $|e_m\rangle\langle e_n|$ to denote the rank one projection $e_k\longmapsto e_m\langle e_n,e_k\rangle$. Let $p:=|e_0\rangle\langle e_0|$ and consider the following operators acting on $\clh$,
\begin{IEEEeqnarray*}{rCl}
a_0=V\otimes 1\otimes 1\quad,\quad b_0=p\otimes U\otimes 1\quad,\quad D_{\theta}=1\otimes e^{-2\sqrt{-1}\pi\theta N} \otimes U\,\,.
\end{IEEEeqnarray*}
Let $C(U_{0,\theta})$ be the $C^*$-subalgebra of $\mathcal{B}(\clh)$ generated by $a_0$, $b_0$ and $D_{\theta}$.

\bppsn\label{isomorphism}
One has $C(U_q(2))=C(U_{0,\theta})$ as $C^*$-algebras.
\eppsn
\prf Observe that $b_0b_0^*=p\otimes 1\times 1$ and thus, $a_0^jb_0b_0^* (a_0^*)^i=p_{ij}\otimes 1\otimes 1$ where $p_{ij}$ is the rank one projection $|e_j\rangle\langle e_i|$ on $\ell^2(\mathbb{N})$. This shows that $\clk\otimes 1\otimes 1\subseteq C(U_{0,\theta})$ where $\clk$ is the space of all compact operators on $\ell^2(\mathbb{N})$. Since 
\begin{IEEEeqnarray}{rCl}\label{used here itself}
a-a_0 &=& \Big(\sqrt{1-|q|^{2N}}-1\Big)V\otimes 1\otimes 1\in\clk\otimes 1\otimes 1\,,
\end{IEEEeqnarray}
we have $a \in C(U_{0,\theta})$. Now, observe that
\begin{IEEEeqnarray*}{rCl}
b &=& \sum_{n=0}^{\infty}(q^n|e_n\rangle\langle e_n|)\otimes U\otimes 1\\
&=& \sum_{n=0}^{\infty}q^na_0^nb_0(a_0^*)^n,
\end{IEEEeqnarray*}
and this is an element in $C(U_{0,\theta})$. Thus, we have $C(U_q(2))\subseteq C(U_{0,\theta})$. To see the reverse inclusion, first observe that 
\[
\mathds{1}_{\{1\}}(b^*b)=p \otimes 1 \otimes 1 \in C(U_q(2)).
\]
Since the following element
\[
a^j(p \otimes 1 \otimes 1)(a^*)^i=\prod_{k=0}^{i-1}\sqrt{1-|q|^{2(i-k)}}\,p_{ij}\otimes 1\otimes 1
\]
is in $C(U_q(2))$, we get that $\clk\otimes 1\otimes 1\subseteq C(U_q(2))$. Since $a-a_0\in\clk\otimes 1\otimes 1$ by Eqn. \ref{used here itself}, we get that $a_0\in C(U_q(2))$. Finally,
\[
 b_0=(p \otimes 1 \otimes 1)(q^N \otimes U\otimes 1)=(p \otimes 1 \otimes 1)b \in C(U_q(2)).
\]
Thus, $C(U_{0,\theta})\subseteq C(U_q(2))$ and this completes the proof.\qed

Let $\scrt:=C^*(V)$ be the Toeplitz algebra. We have the well-known short exact sequence
\[
0 \longrightarrow\mathcal{K}\stackrel{\iota}{\longrightarrow}\scrt\stackrel{\sigma}{\longrightarrow} C(\mathbb{T}) \longrightarrow 0
\]
where $\sigma:V\longmapsto\mathbf{z}$ (here $\mathbf{z}$ denotes the standard unitary generator for $C(\mathbb{T})$). Since $C(U_q(2))$ is a subalgebra of $\scrt\otimes\clb(\ell^2(\bbz)\otimes\ell^2(\bbz))$, we consider the homomorphism $\tau:C(U_q(2))\longrightarrow C(\bbbt)\otimes\clb(\ell^2(\bbz)\otimes\ell^2(\bbz))$ given by $\tau=\sigma\otimes 1\otimes 1$, and let
\begin{IEEEeqnarray*}{rCl}
\mathcal{I}_\theta &=& \mbox{ the closed two-sided ideal of } C(U_q(2)) \mbox{ generated by } b_0 \mbox{ and } b_0^*\,,\\
\mathcal{B}_\theta &=& C^*\big(\{\tau(a_0)\,,\,\tau(D_\theta)\}\big)=C^*\big(\{\mathbf{z}\otimes 1\otimes 1\,,\,1\otimes e^{-2\sqrt{-1}\pi\theta N}\otimes U\}\big)\,.
\end{IEEEeqnarray*}

\bppsn\label{short exact}
The following chain of $C^*$-algebras
\[
0 \longrightarrow\mathcal{I}_\theta\stackrel{\iota}{\longrightarrow} C(U_q(2))\stackrel{\tau}{\longrightarrow}\mathcal{B}_\theta\longrightarrow 0
\]
is an exact sequence, where `$\iota$' denotes the inclusion map.
\eppsn
\prf Since $\tau(b_0)=\tau(b_0^*)=0$, we have $\mathcal{I}_\theta\subseteq\mbox{ ker}(\tau)$. Consider any irreducible representation $\pi$ of $C(U_q(2))$ such that $\mathcal{I}_\theta\subseteq\mbox{ker}(\pi)$, i,e. $\pi(b_0)=0$. By Thm. $3.2$ in \cite{ZhaZha-2005aa}, it follows that $\pi$ is one dimensional. Define $\widetilde{\pi}:\mathcal{B}_\theta\longrightarrow\mathbb{C}$ by $\widetilde{\pi}(\tau(a_0))=\pi(a_0)$ and $\widetilde{\pi}(\tau(D_\theta))=\pi(D_\theta)$. Then, $\pi$ factors through the map $\tau:C(U_q(2))\longrightarrow\mathcal{B}_\theta$. This shows that $\mathcal{I}_\theta=\mbox{ker}(\tau)$.\qed

Let $\mathcal{C}_\theta$ be the $C^*$-subalgebra of $\mathcal{B}(\ell^2(\bbz)\otimes\ell^2(\bbz))$ generated by $U\otimes 1$ and $e^{-2\sqrt{-1}\pi\theta N}\otimes U$. Since $\theta$ is irrational, by the universality and simpleness of 
the noncommutative torus $\mathbb{A}_\theta$, we get that $\mathcal{C}_\theta\cong\mathbb{A}_\theta$ as $C^*$-algebras.
\blmma\label{kernel is NC torus}
$\mathcal{I}_\theta=\mathcal{K}(\ell^2(\bbn))\otimes\mathcal{C}_\theta\,$.
\elmma
\prf We first claim that $\mathcal{K}(\ell^2(\bbn))\otimes \mathcal{C}_\theta\subseteq C(U_q(2))$, and it is a closed two sided ideal. For this, observe that $p\otimes\mathcal{C}_\theta\subseteq C(U_q(2))$ and hence, $a_0^j(p\otimes\mathcal{C}_\theta)a_0^i\subseteq C(U_q(2))$. This is same as $p_{ij}\otimes\mathcal{C}_\theta$ and hence, $\mathcal{K}(\ell^2(\bbn))\otimes\mathcal{C}_\theta\subseteq C(U_q(2))$. It is now a straightforward verification that $\mathcal{K}(\ell^2(\bbn))\otimes\mathcal{C}_\theta$ becomes a two sided ideal in $C(U_q(2))$.

To complete the proof, observe that $b\in\mathcal{K}(\ell^2(\bbn))\otimes\mathcal{C}_\theta$ and hence, $\mathcal{I}\subseteq\mathcal{K}(\ell^2(\bbn))\otimes\mathcal{C}_\theta$. Moreover, the image of $\mathcal{K}(\ell^2(\bbn))\otimes\mathcal{C}_\theta$ under the map $\tau$ is zero and hence, it follows that $\mathcal{K}(\ell^2(\bbn))\otimes\mathcal{C}_\theta\subseteq\mbox{ ker}(\tau)=\mathcal{I}_\theta$ by Propn. \ref{short exact}.\qed

\blmma\label{needed K groups}
For $j=0,1,\,K_j(\mathcal{I}_\theta)$ and $K_j(\mathcal{B}_\theta)$ all are isomorphic to $\bbz^2$. The generators of $K_0(\mathcal{I}_\theta)$ are $[p\otimes 1\otimes 1]$ and $[p\otimes p_\theta]$, where $p_\theta$ is the Powers-Rieffel projection in $\mathbb{A}_\theta$ with trace $\theta$, and that of $K_1(\mathcal{I}_\theta)$ are $[p\otimes U\otimes 1+(1-p)\otimes 1\otimes 1]$ and $[p\otimes e^{-2\sqrt{-1}\pi\theta N}\otimes U+(1-p)\otimes 1\otimes 1]$.
\elmma
\prf Since by the definition $\mathcal{B}_\theta=C^*(\{\tau(a_0),\tau(D_\theta)\})$, it is immediate that $\mathcal{B}_\theta\cong C(\mathbb{T}^2)$. Thus, $K_0(\mathcal{B}_\theta)\cong K_1(\mathcal{B}_\theta)\cong \bbz^2$ and $K_1(\mathcal{B}_\theta)$ is generated by $[\tau(a_0)]$ and $[\tau(D_\theta)]$. Moreover, $K_0(\mathcal{B}_\theta)$ is generated by $[1]$ and $[\mathbb{P}(\tau(a_0),\tau(D_\theta))]$. We refer the reader to  \cite{Con-1981aa} (or Sec. $4$ in \cite{CS-2011aa}) for detail exposition of the operation $\mathbb{P}$. 

By the previous Propn. \ref{kernel is NC torus} and the fact that $\mathcal{C}_\theta$ is isomorphic to the noncommutative torus $\mathbb{A}_\theta$, we get that $K_0(\mathcal{I}_\theta)\cong K_1(\mathcal{I}_\theta)\cong \bbz^2$. Let $p_\theta$ be the Powers-Rieffel projection in the noncommutative torus $\mathbb{A}_\theta$ which generates $K_0(\mathbb{A}_\theta)$ along with $[1]$. Then, the generators of $K_0(\mathcal{I}_\theta)$ are $[p\otimes 1\otimes 1]$ and $[p\otimes p_\theta]$. Moreover, the generators of $K_1(\mathcal{I}_\theta)$ are $[p\otimes U\otimes 1+(1-p)\otimes 1\otimes 1]$ and $[p\otimes e^{-2\sqrt{-1}\pi\theta N}\otimes U+(1-p)\otimes 1\otimes 1]$.\qed

\bthm\label{K groups}
For $q=|q|e^{\sqrt{-1}\pi\theta}$ with $\theta$ irrational, both the $K$-groups $K_{0}(C(U_q(2)))$ and $K_{1}(C(U_q(2)))$ are 
isomorphic to $\bbz^2$. The equivalence classes of unitary $[D]$ and $[p\otimes U\otimes 1+(1-p)\otimes 1\otimes 1]$ form a $\bbz$-basis for $K_{1}(C(U_q(2)))$. The equivalence classes of projections $[1]$ and $[p\otimes p_\theta]$ form a $\bbz$-basis for $K_{0}(C(U_q(2)))$, where $p_\theta$ denotes the Powers-Rieffel projection with trace $\theta$ in the noncommutative torus $\mathbb{A}_\theta$. 
\ethm
\prf In the following six-term exact sequence
\[
\begin{tikzcd}
K_0(\mathcal{I}_\theta) \ar{r}{K_0(\iota)} & K_0(C(U_q(2))) \ar{r}{K_0(\tau)} & K_0(\mathcal{B}_\theta) \ar{d}{\delta} \\
K_1(\mathcal{B}_\theta) \ar{u}{\partial} & K_1(C(U_q(2))) \ar{l}{K_1(\tau)} & K_1(\mathcal{I}_\theta) \ar{l}{K_1(\iota)}
\end{tikzcd}
\]
we first consider the index map $\partial:K_1(\mathcal{B}_\theta)\longrightarrow K_0(\mathcal{I}_\theta)$. Since we 
have the inclusion map in Propn. \ref{short exact}, $\mathcal{I}_\theta$ is an ideal in $C(U_q(2)),\,a_0$ is an isometry 
in $C(U_q(2))$, and $D_\theta$ is unitary in $C(U_q(2))$, we immediately observe that $\partial([\tau(a_0)])=-[p\otimes 1\otimes 1]$ and $\partial([\tau(D_\theta)])=0$ (Propn. $9.2.4$ in \cite{RLL}). So we get that $\mbox{ker}(K_0(\iota))=Im(\partial)=\langle[p\otimes 1\otimes 1]\rangle\cong\mathbb{Z}$ and $Im(K_1(\tau))=\mbox{ker}(\partial)=\langle[\tau(D_\theta)]\rangle\cong\mathbb{Z}$. Hence, $Im(K_0(\iota))=\mbox{ker}(K_0(\tau))\cong\bbz$ by
Lemma \ref{needed K groups}. Now, consider the exponential map $\delta:K_0(\mathcal{B}_\theta)\longrightarrow K_1(\mathcal{I}_\theta)$. Observe that $\delta([1])=0$ as $\tau$ is unital (Propn. $12.2.2$ in \cite{RLL}), and the subgroup generated by the map $\delta([\mathbb{P}(\tau(a_0),\tau(D_\theta))])$ is the subgroup in $K_1(\mathcal{I}_\theta)$ generated by $[D_\theta(1-a_0a_0^*)+a_0a_0^*]$ by Cor. $4.2$ in \cite{CS-2011aa}. But, $D_\theta(1-a_0a_0^*)+a_0a_0^*=p\otimes e^{-2\sqrt{-1}\pi\theta N}\otimes U+(1-p)\otimes 1\otimes 1$. Hence, $\mbox{ker}(\delta)=Im(K_0(\tau))=\langle[1]\rangle\cong\bbz$ and $Im(\delta)=\mbox{ker}(K_1(\iota))=\langle[p\otimes e^{-2\sqrt{-1}\pi\theta N}\otimes U+(1-p)\otimes 1\otimes 1]\rangle\cong\bbz$. Thus, $Im(K_1(\iota))=\mbox{ker}(K_1(\tau))\cong\bbz$ by Lemma \ref{needed K groups}. These 
facts show that $K_{0}(C(U_q(2)))\cong K_{1}(C(U_q(2)))\cong\bbz^2$. The generators of the $K$-groups are obtained using Lemma \ref{needed K groups}. Observe that $K_0(C(U_q(2)))$ will be generated by the generator of $Im(K_0(\iota))\cong K_0(\mathcal{I}_\theta)/Im(\partial)$ and $\mbox{ker}(\delta)$, which are $[p\otimes p_\theta]$ and $[1]$ respectively. Similarly, $K_1(C(U_q(2)))$ will be generated by the generator of $Im(K_1(\iota))\cong K_1(\mathcal{I}_\theta)/Im(\delta)$ and $\mbox{ker}(\partial)$, which are $[p\otimes U\otimes 1+(1-p)\otimes 1\otimes 1]$ and $[D]$ respectively.\qed

Finally, we mention the case when $\theta=0,1$. The deformation parameter $q$ is real in these two situations and $q\in(-1,1)$. 
The $C^*$-algebra $C(U_{0,\theta})$ and Propns. (\ref{isomorphism},\ref{short exact}) make sense and hold in these two cases.  
But the $C^*$-algebra $\,\mathcal{C}_\theta$ in Lemma \ref{kernel is NC torus} becomes $C(\mathbb{T}^2)$ when $\,\theta\in\{0,1\}$. Thus, the generators of $K_0(\mathcal{I}_\theta)$ becomes $[1]$ and $[p\otimes\,Bott]$ where, ``$Bott$'' denotes the Bott projection in $M_2(C(\mathbb{T}^2))$. The proof of Thm. \ref{K groups} holds along the same line, and we get 
that the generators of $K_0(C(U_q(2)))$ are $[1]$ and $[p\otimes\,Bott]$.

%%%%%%%%%%%%%%%%%%%%%%%%%%%%%%%%%%%%%%%%%%%%
%%%%%%%%%%   Action of generators   %%%%%%%%%%%%%%%%%%%%%%%
%%%%%%%%%%%%%%%%%%%%%%%%%%%%%%%%%%%%%%%%%%%%

\newsection{Action of the generators on \texorpdfstring{$L^2(h)$}{}}\label{Sec 4}

In this section we describe the action of the generators $a,a^*,b,b^*,D,D^*$ of $U_q(2)$ on  
the orthonormal basis of $L^2(h)$ consisting of matrix coefficients.  
Let $\pi_h: C(U_q(2))\longrightarrow \clb(L^2(h))$ be the GNS representation associated  
to the Haar state $h$. Since the Haar state on $U_q(2)$ is faithful (Thm. $2.8$ in \cite{GuiSau-2020aa}), the representation $\pi_h$ is faithful. We omit the representation symbol $\pi_h$ and write $x$ instead of $\pi_h(x)$ for the action of $x\in C(U_q(2))$. Recall that $\mathscr{A}_q$ is the $\star$-subalgebra of $C(U_q(2))$ generated by the matrix coefficients of all finite dimensional irreducible representations of $U_q(2)$. 

The Peter-Weyl decomposition obtained in Thm. $4.17$ in \cite{GuiSau-2020aa} says that the following set
\[
\big\{|q|^{-i}\sqrt{|2\ell+1|_{|q|}}\,t^\ell_{i,j}D^k:\ell\in\frac{1}{2}\bbn,\,k\in\bbz\big\}
\]
is an orthonormal basis of $L^2(h)$ where,
\begin{IEEEeqnarray*}{rCl}
t_{i,j}^\ell & = & \sum_{\substack{m+n=\ell-i\\ 0\leq m \leq \ell-j\\ 0\leq n \leq \ell+j}}q^{n(\ell-j-m)}\frac{{2\ell\choose \ell+j}_{|q|^2}^{1/2}}{{2\ell\choose \ell+i}_{|q|^2}^{1/2}}{\ell-j\choose m}_{|q|^2}{\ell+j\choose n}_{|q|^2}a^mc^{\ell-j-m}b^nd^{\ell+j-n}
\end{IEEEeqnarray*}
with $\,c=-\bar{q}Db^*$ and $d=Da^*$. Throughout the article we reserve the following notation,
\begin{IEEEeqnarray*}{rCl}
e^\ell_{i,j,k} &:=& |q|^{-i}\sqrt{|2\ell+1|_{|q|}}\,t^\ell_{i,j}(D^*)^k\\
&=& |q|^{-i}\sqrt{|2\ell+1|_{|q|}}\,t^\ell_{i,j}D^{-k}
\end{IEEEeqnarray*}
for the orthonormal basis of $L^2(h)$. For reader's convenience, we recall the following $q$-Jacobi polynomial expressions of $\,t_{i,j}^\ell D^k$ from Sec. $4.3$ in \cite{GuiSau-2020aa}.

\bthm[\cite{GuiSau-2020aa}]\label{Jacobi}
The matrix coefficients $t_{i,j}^\ell D^{-k}$ are expressed in terms of the little q-Jacobi polynomials in the following way~:
\begin{enumerate}[$(i)$]
\item For the case of $\,i+j\leq 0$ and $\,i\geq j$,
\begin{center}
$a^{-(i+j)}c^{i-j}(\bar{q})^{(j-i)(\ell+j)}\frac{{2\ell\choose \ell+j}^{1/2}_{|q|^2}}{{2\ell\choose \ell+i}^{1/2}_{|q|^2}}{\ell-j\choose i-j}_{|q|^2}\mathcal{P}_{\ell+j}^{(i-j,-i-j)}(bb^*;|q|^2)D^{\ell+j-k}\,;$
\end{center}
\item For the case of $\,i+j\leq 0$ and $\,i\leq j$,
\begin{center}
$a^{-(i+j)}b^{j-i}q^{(i-j)(\ell+i)}\frac{{2\ell\choose \ell+j}^{1/2}_{|q|^2}}{{2\ell\choose \ell+i}^{1/2}_{|q|^2}}{\ell+j\choose j-i}_{|q|^2}\mathcal{P}_{\ell
+i}^{(j-i,-i-j)}(bb^*;|q|^2)D^{\ell+i-k}\,;$
\end{center}
\item For the case of $\,i+j\geq 0$ and $\,i\leq j$,
\begin{center}
$q^{(i-j)(\ell+i)}\frac{{2\ell\choose \ell+j}^{1/2}_{|q|^2}}{{2\ell\choose \ell+i}^{1/2}_{|q|^2}}{\ell+j\choose j-i}_{|q|^2}\mathcal{P}_{\ell-j}^{(j-i,i+j)}(bb^*;|q|^2)(a^*)^{i+j}b^{j-i}D^{\ell+i-k}\,;$
\end{center}
\item For the case of $\,i+j\geq 0$ and $\,i\geq j$,
\begin{center}
$(\bar{q})^{(j-i)(\ell+j)}\frac{{2\ell\choose \ell+j}^{1/2}_{|q|^2}}{{2\ell\choose \ell+i}^{1/2}_{|q|^2}}{\ell-j\choose i-j}_{|q|^2}\mathcal{P}_{\ell-i}^{(i-j,i+j)}(bb^*;|q|^2)(a^*)^{i+j}c^{i-j}D^{\ell+j-k}\,.$
\end{center}
\end{enumerate}
\ethm

For $n,k\in\bbz$ and $m,r\in\bbn$, define
\begin{IEEEeqnarray*}{rCl}
a_nb^m(b^*)^rD^k=\begin{cases}
                            a^n b^m(b^*)^r D^k & \mbox{ if } n \geq 0, \cr
                            (a^*)^{-n} b^m(b^*)^r D^k & \mbox{ if } n < 0. \cr
                           \end{cases}
\end{IEEEeqnarray*}

\bthm[\cite{ZhaZha-2005aa}]\label{basis}
The set $\{a_nb^m(b^*)^rD^k:n,k\in\bbz,\,m,r\in\bbn\}$ forms a linear basis of $\,\mathscr{A}_q$ for all $\,q\in\mathbb{C}^*$.
\ethm

Thm. \ref{Jacobi} and \ref{basis} lead us to the following very important theorem, which is the backbone to our search for Dirac operator in Sec. \ref{Sec 5}.
\bthm\label{the action}
The action of the generators of $U_q(2)$ on the orthonormal basis element $e_{i,j,k}^\ell$ is described by the following~$:$
\begin{enumerate}[$(i)$]
\item $D\,\triangleright\,e_{i,j,k}^\ell=\Big(\frac{q}{\overline{q}}\Big)^{(i-j)}e_{i,j,k-1}^\ell$
\item $D^*\,\triangleright\,e_{i,j,k}^\ell=\Big(\frac{\overline{q}}{q}\Big)^{(i-j)}e_{i,j,k+1}^\ell$
\item $b\,\triangleright\,e_{i,j,k}^\ell=\beta_+(\ell,i,j)e_{i-1/2\,,\,j+1/2\,,\,k}^{\ell+1/2}+\beta_-(\ell,i,j)e_{i-1/2\,,\,j+1/2\,,\,k-1}^{\ell-1/2}$
\item $b^*\,\triangleright\,e_{i,j,k}^\ell=\beta_+^+(\ell,i,j)e_{i+1/2\,,\,j-1/2\,,\,k+1}^{\ell+1/2}+\beta_-^+(\ell,i,j)e_{i+1/2\,,\,j-1/2\,,\,k}^{\ell-1/2}$
\item $a\,\triangleright\,e_{i,j,k}^\ell=\alpha_+(\ell,i,j)e_{i-1/2\,,\,j-1/2\,,\,k}^{\ell+1/2}+\alpha_-(\ell,i,j)e_{i-1/2\,,\,j-1/2\,,\,k-1}^{\ell-1/2}$
\item $a^*\,\triangleright\,e_{i,j,k}^\ell=\alpha_+^+(\ell,i,j)e_{i+1/2\,,\,j+1/2\,,\,k+1}^{\ell+1/2}+\alpha_-^+(\ell,i,j)e_{i+1/2\,,\,j+1/2\,,\,k}^{\ell-1/2}$
\end{enumerate}
where,
\begin{IEEEeqnarray*}{rCl}
\beta_+(\ell,i,j) & = & q^{\ell-j}\sqrt{\frac{(1-|q|^{2(\ell+j+1)})(1-|q|^{2(\ell-i+1)})}{(1-|q|^{2(2\ell+1)})(1-|q|^{2(2\ell+2)})}}\,.
\end{IEEEeqnarray*}
\begin{IEEEeqnarray*}{rCl}
\beta_-(\ell,i,j) &=& -q^{\ell-j-1}(\overline{q})^{j-i+1}\sqrt{\frac{q}{\bar{q}}}\sqrt{\frac{(1-|q|^{2(\ell-j)})(1-|q|^{2(\ell+i)})}{(1-|q|^{2(2\ell)})(1-|q|^{2(2\ell+1)})}}\,.
\end{IEEEeqnarray*}
\begin{IEEEeqnarray*}{rCl}
\beta_+^+(\ell,i,j) &=&  -q^{j-i-1}(\overline{q})^{\ell-j+1}\sqrt{\frac{q}{\bar{q}}}\sqrt{\frac{(1-|q|^{2(\ell-j+1)})(1-|q|^{2(\ell+i+1)})}{(1-|q|^{2(2\ell+1)})(1-|q|^{2(2\ell+2)})}}\,.
\end{IEEEeqnarray*}
\begin{IEEEeqnarray*}{rCl}
\beta_-^+(\ell,i,j) &=&  (\overline{q})^{\ell-j}\sqrt{\frac{(1-|q|^{2(\ell+j)})(1-|q|^{2(\ell-i)})}{(1-|q|^{2(2\ell)})(1-|q|^{2(2\ell+1)})}}\,.
\end{IEEEeqnarray*}
\begin{IEEEeqnarray*}{rCl}
\alpha_+(\ell,i,j) &=& \sqrt{\frac{(1-|q|^{2(\ell-j+1)})(1-|q|^{2(\ell-i+1)})}{(1-|q|^{2(2\ell+1)})(1-|q|^{2(2\ell+2)})}}\,.
\end{IEEEeqnarray*}
\begin{IEEEeqnarray*}{rCl}
\alpha_-(\ell,i,j) &=& q^{\ell-i}(\overline{q})^{\ell-j+1}\sqrt{\frac{q}{\bar{q}}}\sqrt{\frac{(1-|q|^{2(\ell+j)})(1-|q|^{2(\ell+i)})}{(1-|q|^{2(2\ell)})(1-|q|^{2(2\ell+1)})}}\,.
\end{IEEEeqnarray*}
\begin{IEEEeqnarray*}{rCl}
\alpha_+^+(\ell,i,j) &=& q^{\ell-j}(\overline{q})^{\ell-i+1}\sqrt{\frac{q}{\bar{q}}}\sqrt{\frac{(1-|q|^{2(\ell+j+1)})(1-|q|^{2(\ell+i+1)})}{(1-|q|^{2(2\ell+1)})(1-|q|^{2(2\ell+2)})}}\,.
\end{IEEEeqnarray*}
\begin{IEEEeqnarray*}{rCl}
\alpha_-^+(\ell,i,j) &=& \sqrt{\frac{(1-|q|^{2(\ell-j)})(1-|q|^{2(\ell-i)})}{(1-|q|^{2(2\ell)})(1-|q|^{2(2\ell+1)})}}\,.
\end{IEEEeqnarray*}
\ethm
\prf We omit the proof as it follows from Thm. (\ref{Jacobi}, \ref{basis}) by straightforward but tedious algebraic computations.\qed

\brmrk
These coefficients will play a pivotal role in finding equivariant Dirac operator in Sec. \ref{Sec 6}.
\ermrk

\hspace*{-.6cm}\textbf{Notation~:} $1-t^{r\pm s}:=(1-t^{r+s})(1-t^{r-s})$ for $t\in\mathbb{R}_+$.
\bcrlre\label{action of bb*}
We have,
\begin{center}
$bb^*\,\triangleright\,e_{i,j,k}^\ell=\gamma_+(\ell,i,j)e_{i,j\,,\,k+1}^{\ell+1}+\gamma(\ell,i,j)e_{i,j,k}^\ell+\gamma_-(\ell,i,j)e_{i,j\,,\,k-1}^{\ell-1}$
\end{center}
where,
\begin{IEEEeqnarray*}{rCl}
\gamma_+(\ell,i,j) &=& -q^{\ell-i}(\bar{q})^{\ell-j+1}\sqrt{\frac{q}{\bar{q}}}\,\frac{1}{\big(1-|q|^{2(2\ell+2)}\big)}\sqrt{\frac{\big(1-|q|^{2(\ell+1\pm j)}\big)\big(1-|q|^{2(\ell+1\pm i)}\big)}{\big(1-|q|^{2(2\ell+1)}\big)\big(1-|q|^{2(2\ell+3)}\big)}}\\
\gamma(\ell,i,j) &=& |q|^{2(\ell-i)}\frac{\big(1-|q|^{2(\ell-j+1)}\big)\big(1-|q|^{2(\ell+i+1)}\big)}{\big(1-|q|^{2(2\ell+1)}\big)\big(1-|q|^{2(2\ell+2)}\big)}+|q|^{2(\ell-j)}\frac{\big(1-|q|^{2(\ell+j)}\big)\big(1-|q|^{2(\ell-i)}\big)}{\big(1-|q|^{2.2\ell}\big)\big(1-|q|^{2(2\ell+1)}\big)}\\
\gamma_-(\ell,i,j) &=& -q^{\ell-j-1}(\bar{q})^{\ell-i}\sqrt{\frac{q}{\bar{q}}}\,\frac{1}{(1-|q|^{2.2\ell})}\sqrt{\frac{\big(1-|q|^{2(\ell\pm j)}\big)\big(1-|q|^{2(\ell\pm i)}\big)}{\big(1-|q|^{2(2\ell-1)}\big)\big(1-|q|^{2(2\ell+1)}\big)}}\,\,.
\end{IEEEeqnarray*}
\ecrlre

%%%%%%%%%%%%%%%%%%%%%%%%%%%%%%%%%%%%%%%
%%%%%%%%%%   Fixed points   %%%%%%%%%%%%%%%%%%%%%%
%%%%%%%%%%%%%%%%%%%%%%%%%%%%%%%%%%%%%%%

\newsection{Fixed points for the action of \texorpdfstring{$\,bb^*$}{}}\label{Sec 5}

The aim of this section is to analyse the fixed point space for the action of $bb^*$, and this will be needed to compute the index pairing in Sec. \ref{Sec 7}.

It is easy to see that the spectrum of $bb^*$, denoted by $\sigma(bb^*)$, is $\{1,|q|^2,|q|^4,\cdots\} \cup \{0\}$.
Fix $i,j \in \frac{1}{2}\bbz$ and $k \in \bbz$. We abbreviate ``closed linear span'' by c.l.s. For any $\lambda\in\mathbb{C}\setminus\{0\}$, define the followings,
\begin{IEEEeqnarray*}{rCl}
w_{ij} &=& \max \{|i|,|j|\}\\
A(i,j,k) &=& \mbox{c.l.s.}\,\big{\{}e_{i,j,k+m}^{w_{ij}+m}:m\in\bbn\big{\}}\\
E_{\lambda}(i,j,k) &=&\{v \in A(i,j,k): bb^*v=\lambda v\}\\
E_{\lambda} &=& \{ v \in L^2(h): bb^*v=\lambda v\}.
\end{IEEEeqnarray*}
It follows from Cor. \ref{action of bb*} that $ A(i,j,k)$ is an invariant subspace of $bb^*$. Using Thm. $4.17$ in \cite{GuiSau-2020aa}, we get that 
\[
 L^2(h)=\bigoplus A(i,j,k)\,.
\]
Note that either both $i$ and $j$ are integers or both are half integers. 
To get solutions of the equation $bb^*v=\lambda v$, it is enough to look at  solutions in each $ A(i,j,k)$.
Suppose that $v =\sum_{m=0}^{\infty}c_me_{i,j,k+m}^{w_{ij}+m} \in A(i,j,k)$ is a nonzero solution of $bb^*v=\lambda v$. 
By Cor. \ref{action of bb*}, we get the following recurrence relations,
\begin{IEEEeqnarray}{lCl}\label{recurrence}
\lambda c_0 &=& c_0\gamma(w_{ij},i,j)+c_1\gamma_-(w_{ij}+1,i,j)\,,
\end{IEEEeqnarray}
and for $m\geq 1$,
\begin{IEEEeqnarray}{lCl}\label{recurrence1}
\lambda c_m &=& c_{m+1}\gamma_-(w_{ij}+m+1,i,j)+c_m\gamma(w_{ij}+m,i,j)+c_{m-1}\gamma_+(w_{ij}+m-1,i,j)\,. 
\end{IEEEeqnarray}
Observe that in Cor. \ref{action of bb*}, $\gamma_+(\ell,i,j)$ and $\gamma(\ell,i,j)$ are always nonzero, and $\gamma_-(\ell,i,j) \neq 0$ for $w_{ij}<\ell$. Using this, from the relations \ref{recurrence} and \ref{recurrence1} it follows that $c_0 \neq 0$ as otherwise $v$ would be zero. We call $e_{i,j,k}^{w_{ij}}$ the leading term and the corresponding coefficient the leading coefficient.

\blmma\label{few observations}
Let $\lambda \in \sigma(bb^*)$. Then one has the followings,
\begin{enumerate}[$(i)$]
\item $\dim E_{\lambda}(i,j,k)\leq 1$.
\item If $v$ is a nonzero vector in $ E_{\lambda}$, then $Dv$ and $D^*v$ are nonzero vectors in  $ E_{\lambda}$.
\item If $v\in E_{\lambda}$, then $av$ and $a^*v$ are in  $ E_{|q|^2\lambda}$ and $ E_{\frac{\lambda}{|q|^2}}$ respectively.
\item If $v$ is a nonzero vector in $E_{\lambda}(i,j,k)$, then $bv$  is a nonzero vector 
in $E_{\lambda}(i-\frac{1}{2},j+\frac{1}{2},k)$ for $i\leq j$, and in $E_{\lambda}(i-\frac{1}{2},j+\frac{1}{2},k-1)$ for $i > j$.
\item If $ v$ is a nonzero vector in $E_{\lambda}(i,j,k)$, then $b^*v$ is a nonzero vector 
in $E_{\lambda}(i+\frac{1}{2},j-\frac{1}{2},k+1)$ for $i\geq j$, and in $E_{\lambda}(i+\frac{1}{2},j-\frac{1}{2},k)$ for $i< j$.
\item For any nonzero vector $v$ in $E_{\lambda}(i,j,k)$, $av$ is a nonzero vector 
in $E_{|q|^2\lambda}(i-\frac{1}{2},j-\frac{1}{2},k-1)$ for $i\geq j$, and in $E_{|q|^2\lambda}(i-\frac{1}{2},j-\frac{1}{2},k)$ for $i< j$.
\item If $\lambda \neq 1$, then for any nonzero vector $v$ in $E_{\lambda}(i,j,k)$, $a^*v$ is a nonzero vector 
in $E_{\frac{\lambda}{|q|^2}}(i+\frac{1}{2},j+\frac{1}{2},k+1)$  for $i\geq j$, and in $E_{\frac{\lambda}{|q|^2}}(i+\frac{1}{2},j+\frac{1}{2},k)$  for $i<j$.
\end{enumerate}
\elmma
\prf 
\begin{enumerate}[$(i)$]
 \item The eqn. (\ref{recurrence}) and $(\ref{recurrence1})$ have a unique solution $(c_m)_{m\geq 0}$ if one fix the leading coefficient $c_0$.  
 Now if $\sum_{m=0}^{\infty} |c_m|^2 < \infty$, then $\dim E_{\lambda}(i,j,k)= 1$ otherwise $\dim E_{\lambda}(i,j,k)= 0$.
 \item Since $D$ is a unitary, $Dv\neq 0$ if $v \neq 0$. 
 Moreover, $Dv\in E_{\lambda}$ for $v \in E_{\lambda}$ follows from the fact that $D$ commutes with $bb^*$. Similar reason for $D^*$.
 \item Follows from the defining relations $ba=qab$ and $a^*b=qba^*$. 
 \item That $v\neq 0$ in $E_{\lambda}(i,j,k)$ implies $bv\neq 0$ follows immediately from the normality of $b$. The rests follow by looking at the leading coefficient of $v$ and the action of $b$ given by Thm. \ref{the action}.
 \item Use the same argument as in part $(iv)$. 
 \item Let $av=0$. Using the relation $a^*a+ |q|^2bb^*=1$, one has $bb^*v=\frac{1}{|q|^2}v$. 
 Since $\frac{1}{|q|^2} \notin \sigma(bb^*)$, one arrives at a contradiction. The rest will follow by analysing the 
 leading coefficient of $v$ and the action of $a^*$ given by Thm. \ref{the action}.
 \item Use the same argument as in part $(vi)$.\qed
\end{enumerate}

Now, we provide a nontrivial solution of the equation $bb^*v=v$. More precisely, we show that $\dim E_{1}(0,0,0)=1$.
\blmma \label{solution_00}
Let $v=\sum_{m=0}^\infty c_me^m_{0,0,m}\in L^2(h)$. The equation $bb^*v=v$ has a unique nonzero solution up to constant scalar multiple.
\elmma
\prf By Cor. \ref{action of bb*}, we see that
\begin{IEEEeqnarray*}{rCl}
bb^*(e^m_{0,0,m}) &=& \sum_{\xi=-1}^1\,\Upsilon_\xi(m)e^{m+\xi}_{0,0,m+\xi}\,
\end{IEEEeqnarray*}
where,
\begin{IEEEeqnarray*}{rCl}
\Upsilon_{-1}(m) &=& -\frac{|q|^{2m-1}(1-|q|^{2m})^2}{(1-|q|^{4m})\sqrt{(1-|q|^{4m-2})(1-|q|^{4m+2})}}\,\,,\\
\Upsilon_0(m) &=& \frac{|q|^{2m}(1-|q|^{2m+2})^2}{(1-|q|^{4m+2})(1-|q|^{4m+4})}+\frac{|q|^{2m}(1-|q|^{2m})^2}{(1-|q|^{4m})(1-|q|^{4m+2})}\,\,,\\
\Upsilon_1(m) &=& -\frac{|q|^{2m+1}(1-|q|^{2m+2})^2}{(1-|q|^{4m+4})\sqrt{(1-|q|^{4m+2})(1-|q|^{4m+6})}}\,\,.
\end{IEEEeqnarray*}
We want to solve $bb^*v=v$ i,e. $\sum_{m=0}^\infty c_me^m_{0,0,m}=\sum_{m=0}^\infty\sum_{\xi=-1}^1c_m\Upsilon_\xi(m) e^{m+\xi}_{0,0,m+\xi}$. That is,
\begin{IEEEeqnarray*}{rCl}
\sum_{m=0}^\infty c_me^m_{0,0,m} &=& \sum_{m=0}^\infty c_m\Upsilon_{-1}(m)e^{m-1}_{0,0,m-1}+c_m\Upsilon_0(m)e^m_{0,0,m}+c_m\Upsilon_1(m)e^{m+1}_{0,0,m+1}\,.
\end{IEEEeqnarray*}
Equating coefficients of the basis elements of both sides, and observing that $\Upsilon_1(m-1)=\Upsilon_{-1}(m)$, we get the following,
\begin{IEEEeqnarray*}{rCl}
c_0 &=& c_0\Upsilon_0(0)+c_1\Upsilon_{-1}(1)\,,\\
c_1 &=& c_2\Upsilon_{-1}(2)+c_1\Upsilon_0(1)+c_0\Upsilon_1(0)\\
&=& c_2\Upsilon_1(1)+c_1\Upsilon_0(1)+c_0\Upsilon_1(0)\,,\\
c_m &=& c_{m+1}\Upsilon_1(m)+c_m\Upsilon_0(m)+c_{m-1}\Upsilon_1(m-1)\quad\forall\,m\geq 2\,.
\end{IEEEeqnarray*}
Observe that if $c_0=0$ then $c_m=0$ for all $m>0$, and consequently $v=0$. Assume that $c_0=1$. Then, $c_1=\frac{1-\Upsilon_0(0)}{\Upsilon_1(0)}$ and we have the following recurrence relation,
\begin{IEEEeqnarray}{rCl}\label{recurrence_0}
c_{m+1} &=& \frac{c_m(1-\Upsilon_0(m))-c_{m-1}\Upsilon_1(m-1)}{\Upsilon_1(m)}\quad\forall\,m\geq 1\,.
\end{IEEEeqnarray}
Plugging the value of $\Upsilon_0(0)$ and $\Upsilon_1(0)$ we see that $c_1=-|q|\sqrt{\frac{1-|q|^6}{1-|q|^2}}\,$. Now, one can verify that for $m\geq 2$,
\[
c_m=(-1)^{m+1}c_1|q|^{m^2-1}\sqrt{\frac{1-|q|^{4m+2}}{1-|q|^6}}
\]
satisfies the recurrence relation in Eqn. \ref{recurrence_0}. That is, for all $m\geq 0$ we have,
\[
c_m=(-1)^m|q|^{m^2}\sqrt{\frac{1-|q|^{4m+2}}{1-|q|^2}}\,.
\]
Now, observe that
\begin{IEEEeqnarray*}{rCl}
\frac{|c_m|}{|q|^m} &=& |q|^{m^2-m}\sqrt{\frac{1-|q|^{4m+2}}{1-|q|^2}}\\
&\longrightarrow& \,\,0\quad\mbox{as}\quad m\rightarrow\infty
\end{IEEEeqnarray*}
because $0<|q|<1$. Thus, $c_m=o(|q|^m)$ as $m\rightarrow\infty$, i,e. $(c_m)_{m\geq0}\in\ell^2(\bbn)$. This completes the proof.\qed

Let $r\in\bbn$ and consider the operator $bb^*$ restricted to the invariant  subspace $A(0,0,0)$. Consider the commutative unital $C^*$-algebra $C^*(1,bb^*)\cong C(\sigma(bb^*))$. By the spectral decomposition, there exists a subset $\Omega=\{n_r:r\in\bbn,\,n_r<n_{r+1}\}$ of $\bbn$ such that $E_{|q|^{2n_r}}(0,0,0) \neq \{0\}$, and we can write 
\[
 A(0,0,0)= E_{0}(0,0,0) \bigoplus_{r \in \bbn}   E_{|q|^{2n_r}}(0,0,0).
\]
By Lemma \ref{solution_00}, we have $n_0=0$. It follows from part $(i)$ of the Lemma \ref{few observations} that the cardinality of $\Omega$ is infinite.  

\bppsn\label{E_1}
One has the followings.
\begin{enumerate}[$(i)$]
\item If $\dim E_1(i,j,k_0)=1$ for some $k_0\in\bbz$, then for each $k\in\bbz$ we have $\dim E_1(i,j,k)=1$.
\item If $\dim E_1(i,j,k)=1$, then $\dim E_1(i^{'},j^{'},k^{'})=1$ for $k^{'}\in \bbz$, and $i^{'},j^{'}\in\frac{1}{2}\bbz$ with $i^{'}+j^{'}=i+j$.
\item For $n_r\in\Omega,\,\dim E_1(i,j,k)=1$ for all $k\in\bbz$ and $i,j\in\frac{1}{2}\bbz$ such that $i+j=n_r$.
\item $E_1=\bigoplus_{r=0}^{\infty}\bigoplus_{i+j=n_r\,,\,k\in\bbz}E_1(i,j,k)$.
\end{enumerate}
\eppsn
\prf
\begin{enumerate}[$(i)$]
\item It is a direct consequence of part $(ii)$ of Lemma \ref{few observations}.
\item This follows by repeatedly applying part $(ii),\,(iv)$ and $(v)$ of Lemma \ref{few observations}, together with part $(i)$.
\item We have $\dim E_{|q|^{2n_r}}(0,0,0)=1$. Using part $(vii)$ of Lemma \ref{few observations} repeatedly, we get that $\dim E_{1}(\frac{n_r}{2},\frac{n_r}{2},n_r)=1$. By part $(ii)$, the claim now follows.
\item From part $(iii)$, we have $E_1\supseteq\bigoplus_{r=0}^{\infty}\bigoplus_{i+j=n_r\,,\,k\in\bbz}E_1(i,j,k)$. To prove the claim, it is enough to show that $\dim E_1(i,j,0)=0$ for $i,j\in\frac{1}{2}\bbz$ such that $i+j<0$. Let $i_0,j_0\in\frac{1}{2}\bbz$ with $i_0+j_0<0$, and let $v=\sum_{m=0}^{\infty}c_me_{i_0,j_0,m}^{w_{i_0j_0}+m}\in A(i_0,j_0,0)$ be a nonzero solution of $bb^*v=v$. By Thm. \ref{the action}, we get the following,
\[
a^*v = c_0 \alpha_-^{+}(w_{i_0j_0},i_0,j_0)e_{i_0+\frac{1}{2},j_0+\frac{1}{2},0}^{w_{i_0j_0}-\frac{1}{2}}+\sum_{m=1}^{\infty}\widetilde{c_m}e_{i_0+\frac{1}{2},j_0+\frac{1}{2},m}^{w_{i_0j_0}-\frac{1}{2}+m}\,,
\]
where $\widetilde{c_m}$'s are  complex numbers. Since $i_0+j_0<0$, we have $i_0\,,\,j_0\neq w_{i_0j_0}$. Hence, $\alpha_-^{+}(w_{i_0j_0},i_0,j_0) \neq 0$, which implies that $a^*v$ is a nonzero vector as $c_0\neq 0$. This along with part $(iii)$ of the Lemma \ref{few observations} gives that $\frac{1}{|q|^2}\in\sigma(bb^*)$, which is a contradiction.\qed
\end{enumerate}

%%%%%%%%%%%%%%%%%%%%%%%%%%%%%%%%%%%%%%%
%%%%%%%%%%   Dirac Operator   %%%%%%%%%%%%%%%%%%%%%
%%%%%%%%%%%%%%%%%%%%%%%%%%%%%%%%%%%%%%%

\newsection{Equivariant Dirac operator}\label{Sec 6}

In this section, we construct a finitely summable Dirac operator on $U_q(2)$ which is equivariant under the comultiplication action, as defined in Section $2$. Let $T$ be the following unbounded operator on $L^2(h)$ with dense domain $\mathscr{A}_q$, defined by
\[
T(e_{i,j,k}^{\ell})=d(\ell,i,k)e_{i,j,k}^{\ell}
\]
where,
\begin{IEEEeqnarray}{lCl}\label{the sequence}
d(\ell,i,k)=\begin{cases}
(2\ell+1)+\sqrt{-1}(k-\ell-i)& \mbox{ if }\,\, i\neq-\ell, \cr
-(2\ell+1)+\sqrt{-1}k& \mbox{ if }\,\, i=-\ell. \cr
\end{cases}
\end{IEEEeqnarray}

\blmma\label{bounded commutator}
For each $x\in\mathscr{A}_q$, the operators $[T,x]$ and $[T^*,x]$ initially defined on $\mathscr{A}_q$ extend to bounded operators acting on $L^2(h)$.
\elmma
\prf It is enough to show that $[T,x]$ where $x\in\{a,a^*,b,b^*,D,D^*\}\subseteq\mathscr{A}_q$ extends to a bounded operator on $L^2(h)$. From Thm. \ref{the action} we have the followings,
\begin{IEEEeqnarray*}{lCl}
[T,D](e_{i,j,k}^{\ell}) &=& \big(d(\ell,i,k-1)-d(\ell,i,k)\big)\Big(\frac{q}{\overline{q}}\Big)^{i-j}e_{i,j,k-1}^{\ell}\,,
\end{IEEEeqnarray*}
\begin{IEEEeqnarray*}{lCl}
[T,D^*](e_{i,j,k}^{\ell}) &=& \big(d(\ell,i,k+1)-d(\ell,i,k)\big)\Big(\frac{\overline{q}}{q}\Big)^{i-j}e_{i,j,k+1}^{\ell}\,,
\end{IEEEeqnarray*}
\begin{IEEEeqnarray*}{lCl}
[T,a](e_{i,j,k}^{\ell}) &=& \big(d(\ell+\frac{1}{2},i-\frac{1}{2},k)-d(\ell,i,k)\big)\alpha_+(\ell,i,j)e_{i-\frac{1}{2},j-\frac{1}{2},k}^{\ell+\frac{1}{2}}\\
&  & +\big(d(\ell-\frac{1}{2},i-\frac{1}{2},k-1)-d(\ell,i,k)\big)\alpha_-(\ell,i,j)e_{i-\frac{1}{2},j-\frac{1}{2},k-1}^{\ell-\frac{1}{2}}\,,
\end{IEEEeqnarray*}
\begin{IEEEeqnarray*}{lCl}
[T,a^*](e_{i,j,k}^{\ell}) &=& \big(d(\ell+\frac{1}{2},i+\frac{1}{2},k+1)-d(\ell,i,k)\big)\alpha_+^+(\ell,i,j)e_{i+\frac{1}{2},j+\frac{1}{2},k+1}^{\ell+\frac{1}{2}}\\
&  & +\big(d(\ell-\frac{1}{2},i+\frac{1}{2},k)-d(\ell,i,k)\big)\alpha_-^+(\ell,i,j)e_{i+\frac{1}{2},j+\frac{1}{2},k}^{\ell-\frac{1}{2}}\,,
\end{IEEEeqnarray*}
\begin{IEEEeqnarray*}{lCl}
[T,b](e_{i,j,k}^{\ell}) &=& \big(d(\ell+\frac{1}{2},i-\frac{1}{2},k)-d(\ell,i,k)\big)\beta_+(\ell,i,j)e_{i-\frac{1}{2},j+\frac{1}{2},k}^{\ell+\frac{1}{2}}\\
&  & +\big(d(\ell-\frac{1}{2},i-\frac{1}{2},k-1)-d(\ell,i,k)\big)\beta_-(\ell,i,j)e_{i-\frac{1}{2},j+\frac{1}{2},k-1}^{\ell-\frac{1}{2}}\,,
\end{IEEEeqnarray*}
\begin{IEEEeqnarray*}{lCl}
[T,b^*](e_{i,j,k}^{\ell}) &=& \big(d(\ell+\frac{1}{2},i+\frac{1}{2},k+1)-d(\ell,i,k)\big)\beta_+^+(\ell,i,j)e_{i+\frac{1}{2},j-\frac{1}{2},k+1}^{\ell+\frac{1}{2}}\\
&  & +\big(d(\ell-\frac{1}{2},i+\frac{1}{2},k)-d(\ell,i,k)\big)\beta_-^+(\ell,i,j)e_{i+\frac{1}{2},j-\frac{1}{2},k}^{\ell-\frac{1}{2}}\,.
\end{IEEEeqnarray*}
As $|d(\ell,i,k\pm 1)-d(\ell,i,k)|=1$, it is easy to see that $[T,D]$ and $[T,D^*]$ extends to bounded operators on $L^2(h)$. Now, consider the right hand side of $[T,a](e_{i,j,k}^{\ell})$. We have from \ref{the sequence}
\[
\big|d\big(\ell+\frac{1}{2},i-\frac{1}{2},k\big)-d(\ell,i,k)\big|=1\,,
\]
\[
\big|d\big(\ell-\frac{1}{2},i-\frac{1}{2},k-1\big)-d(\ell,i,k)\big|=\begin{cases}
4\ell+1 & \mbox{ for }\,\,\ell+i=0,1\,;\\
1 & \mbox{ for }\,\,\ell+i\neq 0,1\,.
\end{cases}
\]
We also have from Thm. \ref{the action}
\begin{IEEEeqnarray*}{lCl}
1-|q|^2\leq\alpha_+(\ell,i,j) &=& \sqrt{\frac{(1-|q|^{2(\ell-j+1)})(1-|q|^{2(\ell-i+1)})}{(1-|q|^{2(2\ell+1)})(1-|q|^{2(2\ell+2)})}}\\
&\leq& \frac{(1-|q|^{2(2\ell+1)})}{\sqrt{(1-|q|^{2(2\ell+1)})(1-|q|^{2(2\ell+2)})}}\\
&\leq& 1\,;
\end{IEEEeqnarray*}
and
\begin{IEEEeqnarray*}{lCl}
max_j\,|\alpha_-(\ell,i,j)| &=& |q|^{\ell-i+1}\sqrt{\frac{(1-|q|^{2(\ell+i)})}{(1-|q|^{2(2\ell+1)})}}\,.
\end{IEEEeqnarray*}
Hence,
\begin{IEEEeqnarray*}{lCl}
&  & \big|\big(d\big(\ell-\frac{1}{2},i-\frac{1}{2},k-1\big)-d(\ell,i,k)\big)\alpha_-(\ell,i,j)\big|\\
&\leq& \sqrt{\frac{(1-|q|^{2(\ell+i)})}{(1-|q|^{2(2\ell+1)})}}|q|^{\ell-i}\big|d\big(\ell-\frac{1}{2},i-\frac{1}{2},k-1\big)-d(\ell,i,k)\big|\\
&\leq& \begin{cases}
|q|^{2\ell-1}(4\ell+1) & \mbox{ for }\,\,\ell+i=1;\\
1 & \mbox{ otherwise}.
\end{cases}
\end{IEEEeqnarray*}
This gives the boundedness of the commutator $[T,a]$. The boundedness of $[T,b]$ follows similarly.

Now, consider the right hand side of $[T,b^*](e_{i,j,k}^{\ell})$. We have from \ref{the sequence}
\[
\big|d(\ell-\frac{1}{2},i+\frac{1}{2},k)-d(\ell,i,k)\big|=1\,,
\]
\[
\big|d\big(\ell+\frac{1}{2},i+\frac{1}{2},k+1\big)-d(\ell,i,k)\big|=\begin{cases}
4\ell+3 & \mbox{ for }\,\,\ell+i=0,-1\,;\\
1 & \mbox{ for }\,\,\ell+i\neq 0,-1\,.
\end{cases}
\]
We also have from Thm. \ref{the action}
\begin{IEEEeqnarray*}{lCl}
\beta_-^+(\ell,i,j) &=& |q|^{\ell-j}\sqrt{\frac{(1-|q|^{2(\ell+j)})(1-|q|^{2(\ell-i)})}{(1-|q|^{2(2\ell)})(1-|q|^{2(2\ell+1)})}}\\
&\leq& \sqrt{\frac{1-|q|^{2(\ell-i)}}{1-|q|^{2(2\ell+1)}}}\\
&\leq& 1\,;
\end{IEEEeqnarray*}
and
\begin{IEEEeqnarray*}{lCl}
max_j\,|\beta_+^+(\ell,i,j)| &=& |q|^{\ell-i}\sqrt{\frac{(1-|q|^{2(\ell+i+1)})}{(1-|q|^{2(2\ell+2)})}}\,.
\end{IEEEeqnarray*}
Hence,
\begin{IEEEeqnarray*}{lCl}
&  & \big|\big(d(\ell+\frac{1}{2},i+\frac{1}{2},k+1)-d(\ell,i,k)\big)\beta_+^+(\ell,i,j)\big|\\
&\leq& \sqrt{\frac{(1-|q|^{2(\ell+i+1)})}{(1-|q|^{2(2\ell+2)})}}|q|^{\ell-i}\big|d\big(\ell+\frac{1}{2},i+\frac{1}{2},k+1\big)-d(\ell,i,k)\big|\\
&\leq& \begin{cases}
|q|^{2\ell}(4\ell+3) & \mbox{ for }\,\,\ell+i=0;\\
1 & \mbox{ otherwise}.
\end{cases}
\end{IEEEeqnarray*}
This gives the boundedness of the commutator $[T,b^*]$. The boundedness of $[T,a^*]$ follows similarly.\qed

Define the faithful representation $\pi_{eq}$ of $\,C(U_q(2))$ on $\mathscr{H}:=L^2(h)\otimes\mathbb{C}^2$ by
\[
\pi_{eq}(x)=\left[ {\begin{matrix}
   \pi_h(x)  & 0\\
  0 & \pi_h(x) \\
  \end{matrix} } \right]\,\,.
\]
Let 
\[
\mathscr{D}=\left[ {\begin{matrix}
   0  & T^*\\
  T & 0 \\ 
  \end{matrix} } \right]\quad \mbox{ and }\quad \gamma= \left[ {\begin{matrix}
   1 & 0\\
  0 & -1 \\
  \end{matrix} } \right].
\]
Immediately from Lemma \ref{bounded commutator}, it follows that $[\mathscr{D},\pi_{eq}(x)]$ extends to a bounded operator on $\mathscr{H}$ for each $x\in\mathscr{A}_q$.

\blmma\label{compact resolvent}
$\mathscr{D}$ has compact resolvent.
\elmma
\prf Recall from Eqn. \ref{the sequence} the following,
\begin{IEEEeqnarray*}{lCl}
d(\ell,i,k)=\begin{cases}
(2\ell+1)+\sqrt{-1}(k-\ell-i)& \mbox{ if }\,\, i\neq-\ell, \cr
-(2\ell+1)+\sqrt{-1}k& \mbox{ if }\,\, i=-\ell, \cr
\end{cases}
\end{IEEEeqnarray*}
where $\,\ell\in\frac{1}{2}\bbn$, with $i,j\in\{-\ell,-\ell+1,\ldots,\ell\}$, and $k\in\bbz$. Thus, $|d(\ell,i,k)|^2=(2\ell+1)^2+(k-\ell-i)^2$ regardless of whether $i+\ell=0$ or $i+\ell\neq 0$. For any $k\in\bbz$, we have $|d(\ell,i,k)|^2\geq(2\ell+1)^2$, and consequently $\frac{1}{|d(\ell,i,k)|^2}\leq\frac{1}{(2\ell+1)^2}$ for all $\,\ell\in\frac{1}{2}\bbn$. Now, let $\ell\in\frac{1}{2}\bbn$. We have for any $k\in\bbz$ the following,
\[
0\leq\ell+i\leq 2\ell\quad\mbox{i,e.}\quad k\geq k-(\ell+i)\geq k-2\ell\,.
\]
If $k\geq 2\ell$, then $k-(\ell+i)\geq k-2\ell$ implies that $(k-(\ell+i))^2\geq(k-2\ell)^2$. This says that $\frac{1}{|d(\ell,i,k)|^2}\leq\frac{1}{(2\ell+1)^2+(k-2\ell)^2}$. If $k\leq 0$, then $k\geq k-(\ell+i)$ implies that $(k-(\ell+i))^2\geq k^2$. This says that $\frac{1}{|d(\ell,i,k)|^2}\leq\frac{1}{(2\ell+1)^2+k^2}$. That is,
 \begin{IEEEeqnarray*}{lCl}
\frac{1}{|d(\ell,i,k)|^2}\leq\begin{cases}
\frac{1}{(2\ell+1)^2+(k-2\ell)^2}& \mbox{ if }\,\,k\geq 2\ell, \cr
\frac{1}{(2\ell+1)^2+k^2}& \mbox{ if }\,\,k\leq 0. \cr
\end{cases}
\end{IEEEeqnarray*}
This shows that $|\mathscr{D}|^{-1}$ is a compact operator.\qed

\bthm\label{the spectral triple}
The tuple $(\mathscr{A}_q\,,\mathscr{H},\pi_{eq}\,,\mathscr{D}\,,\gamma)$ is a $4^+$-summable, non-degenerate, even spectral triple on $U_q(2)$ that is equivariant under its own comultiplication action.
\ethm
\prf Combining Lemmas (\ref{bounded commutator},\,\ref{compact resolvent}) it follows that $\mathscr{D}$ is a Dirac operator. Equivariance of $\mathscr{D}$ follows from Propn. \ref{equivariance}. To prove $4^+$-summability, fix $n\in\bbn$ and let
\[
S_n=\big{\{}e_{i,j,k}^{\ell}: -\ell \leq i,j \leq \ell, \ell \in \frac{1}{2}\bbn, k \in\bbz, |d(\ell,i,k)|\leq n\big{\}}. 
\]
From Eqn. \ref{the sequence}, it follows that if $|d(\ell,i,k)|\leq n$, then $\ell \leq n$, $ -n\leq i \leq n$, and hence $-n \leq k \leq 3n$. If one counts the number of such matrix coefficients, then we get the following,
\[
\# S_n\leq 4n\big(1+2^2+ \cdots +(2n+1)^2\big).
\]
To get the lower bound, fix $n\in\bbn$ and consider the following subset of $S_n$,
\[
R_n:=\left\{e_{i,j,k}^{\ell}\in S_n:0\leq i,j\leq\ell;\,0\leq\ell\leq\Big\lfloor\frac{n}{8}\Big\rfloor;\,\Big\lfloor\frac{n}{4}\Big\rfloor+1\leq k\leq\Big\lfloor\frac{n}{2}\Big\rfloor\right\}. 
\]
%Observe that for $e_{i,j,k}^{\ell}\in R_n$, we have $\ell+k+i+1\leq n$ and $k-\ell-i\geq 0$. Therefore,
%\[
% |d(\ell,i,k)|=\sqrt{(2\ell+1)^2+(k-\ell-i)^2}\leq 2\ell+1+k-\ell-i=k+\ell-i+1\leq \ell+k+i+1\leq n;
%\]
%This shows that $R_n \subseteq S_n$ for all $n \in \bbn$, and
Hence, we have
\[
\# S_n\geq \# R_n=\left(\Big\lfloor\frac{n}{4}\Big\rfloor+1\right)\left(1+2^2+\cdots+\Big\lfloor\frac{n}{4}\Big\rfloor^2\right),
\]
which proves the required summability.

Now, we show the non-degeneracy of the spectral triple. Recall that a spectral triple such that $[\mathscr{D},x]=0$ implies $\,x\in\bbc$ is called non-degenerate. In view of Thm. \ref{basis}, consider any finite subset $\,F\subseteq\bbz\times\bbn\times\bbn\times\bbz$ and scalars $\,C_{n,m,r,k}\neq 0$, and take
\[
x=\sum_{(n,m,r,k)\in F}C_{n,m,r,k}\,a_nb^m(b^*)^rD^k
\]
in $\mathscr{A}_q$ such that $x\notin\bbc$. We claim that $[T,x]\neq 0$. Let $\eta=max\{|n|+m+r+|k|:(n,m,r,k)\in F\}$. Then, $\eta>0$ as $x\notin\bbc$. Since $F$ is a finite set, there exists a tuple $(n_0,m_0,r_0,k_0)\in F$ such that $\,\eta=|n_0|+m_0+r_0+|k_0|$. Fix this tuple and we rename it as $(n,m,r,k)$ for notational brevity.\\
\textbf{Case 1~:} $n\geq 0$.\\
Fix any non-zero $\ell\in\frac{1}{2}\bbn$ and take $i=\frac{r-m-n}{2},\,j=\frac{m-r-n}{2}$. It is easy to check that
\[
-\big(\ell+\frac{\eta-|k|}{2}\big)\leq\, i\,,\,j\,\leq \ell+\frac{\eta-|k|}{2}\,.
\]
Now, observe that
\begin{IEEEeqnarray*}{lCl}
&  & \Big\langle e^{\ell+\frac{\eta-|k|}{2}}_{i,j,r-k}\,,\,[T,x]\big(e^\ell_{0,0,0}\big)\Big\rangle\\
&=& \Big\langle e^{\ell+\frac{\eta-|k|}{2}}_{i,j,r-k}\,,\,\big[T,\sum_{(n^\prime,m^\prime,r^\prime,k^\prime)\in F}C_{n^\prime,m^\prime,r^\prime,k^\prime}\,a_{n^\prime}b^{m^\prime}(b^*)^{r^\prime}D^{k^\prime}\big]\big(e^\ell_{0,0,0}\big)\Big\rangle\\
&=& C_{n,m,r,k}\Big\langle e^{\ell+\frac{\eta-|k|}{2}}_{i,j,r-k}\,,\,[T,a_nb^m(b^*)^rD^k]\big(e^\ell_{0,0,0}\big)\Big\rangle\\
&=& C_{n,m,r,k}\prod_{t=1}^r\beta^+_+\Big(\ell+\frac{t-1}{2},\frac{t-1}{2},-\frac{t-1}{2}\Big)\prod_{t=r+1}^{r+m}\beta_+\Big(\ell+\frac{t-1}{2},r-\frac{t-1}{2},-r+\frac{t-1}{2}\Big)\\
&  & \prod_{t=r+m+1}^{\eta-|k|}\alpha_+\Big(\ell+\frac{t-1}{2},\frac{2r+1-t}{2},-\frac{t-2m-1}{2}\Big)\Big(d\Big(\ell+\frac{\eta-|k|}{2},i,r-k\Big)-d(\ell,0,0)\Big)\,.
\end{IEEEeqnarray*}
We claim that this number is nonzero. For this, in view of Thm. \ref{the action}, we only have to check whether the last term
\[
\xi:=d\Big(\ell+\frac{\eta-|k|}{2},i,r-k\Big)-d(\ell,0,0)
\]
is nonzero. Observe that $\ell>0$ by assumption, and hence $i\neq-\ell-\frac{\eta-|k|}{2}$. Thus, from \ref{the sequence} we get that
\[
\xi=\eta-|k|+\sqrt{-1}\Big(r-k-i-\frac{\eta-|k|}{2}\Big)\,.
\]
If $\eta\neq|k|$, then $\Re(\xi)\neq 0$ and we are done. If $\eta=|k|$ then we have $n=m=r=0$ and thus, $i=0$. This gives $\Im(\xi)=-k\neq 0$, as $x\notin\bbc$. So, in both the situations we have $\xi\neq 0$ and consequently, $[T,x]\neq 0$.\\
\textbf{Case 2~:} $n<0$.\\
For the same $\ell,i,j$ chosen in the previous case, observe that
\begin{IEEEeqnarray*}{lCl}
&  & \Big\langle e^{\ell+\frac{\eta-|k|}{2}}_{i,j,r-k-n}\,,\,[T,x]\big(e^\ell_{0,0,0}\big)\Big\rangle\\
&=& \Big\langle e^{\ell+\frac{\eta-|k|}{2}}_{i,j,r-k-n}\,,\,\big[T,\sum_{(n^\prime,m^\prime,r^\prime,k^\prime)\in F}C_{n^\prime,m^\prime,r^\prime,k^\prime}\,a_{n^\prime}b^{m^\prime}(b^*)^{r^\prime}D^{k^\prime}\big]\big(e^\ell_{0,0,0}\big)\Big\rangle\\
&=& C_{n,m,r,k}\Big\langle e^{\ell+\frac{\eta-|k|}{2}}_{i,j,r-k-n}\,,\,[T,a_nb^m(b^*)^rD^k]\big(e^\ell_{0,0,0}\big)\Big\rangle\\
&=& C_{n,m,r,k}\prod_{t=1}^r\beta^+_+\Big(\ell+\frac{t-1}{2},\frac{t-1}{2},-\frac{t-1}{2}\Big)\prod_{t=r+1}^{r+m}\beta_+\Big(\ell+\frac{t-1}{2},r-\frac{t-1}{2},-r+\frac{t-1}{2}\Big)\\
&  & \prod_{t=r+m+1}^{\eta-|k|}\alpha_+^+\Big(\ell+\frac{t-1}{2},-\frac{2m+1-t}{2},\frac{t-2r-1}{2}\Big)\Big(d\Big(\ell+\frac{\eta-|k|}{2},i,r-k-n\Big)-d(\ell,0,0)\Big)\,.
\end{IEEEeqnarray*}
Now, consider 
\[
\xi:=d\Big(\ell+\frac{\eta-|k|}{2},i,r-k-n\Big)-d(\ell,0,0)\,.
\]
From Eqn. \ref{the sequence} we get that
\[
\xi=\eta-|k|+\sqrt{-1}\Big(r-k-n-i-\frac{\eta-|k|}{2}\Big)\,.
\]
Again, if $\eta\neq|k|$, then we are done. If $\eta=|k|$, by similar argument as in the previous case we get that $\Im(\xi)=-k\neq 0$. So, $[T,x]\neq 0$.\qed

%%%%%%%%%%%%%%%%%%%%%%%%%%%%%%%%%%%%%%%
%%%%%%%%%%   Chern Character  %%%%%%%%%%%%%%%%%%%%
%%%%%%%%%%%%%%%%%%%%%%%%%%%%%%%%%%%%%%%

\newsection{Nontriviality of the $\mathcal{K}$-homology class}\label{Sec 7}

The aim of this section is to prove the following theorem.

\bthm\label{final thm}
The Chern character of the spectral triple $(\mathscr{A}_q\,,\mathscr{H},\pi_{eq}\,,\mathscr{D}\,,\gamma)$ is nontrivial.
\ethm

We begin by briefly recalling the standard spectral triple on the noncommutative $2$-torus. We denote by $u$ and $v$ the generating unitaries of $\mathbb{A}_{\theta}$ satisfying the relation $uv=e^{2\pi\sqrt{-1}\theta}vu$. Let $\mathcal{A}$ be the unital $\star$-subalgebra of $\mathbb{A}_{\theta}$ generated by $u$ and $v$. Consider the following representation,
\begin{align*}
\tau:\cla_{\theta} &\longrightarrow \mathcal{B}(\ell^2(\bbz^2))\\
u\longmapsto U\otimes 1\quad &, \quad v\longmapsto e^{-2\pi\sqrt{-1}\theta N}\otimes U\,.
\end{align*}
For a complex valued function $f$ on $\bbz^2$, define the operator $T_f$ acting on $\ell^2(\bbz^2)$ as
\[
T_f(e_m\otimes e_n):=f(m,n)(e_m\otimes e_n)\,.
\]
Then, the following tuple
\[
\left(\mathcal{A},\ell^2(\bbz^2)\otimes\mathbb{C}^2,\left[{\begin{smallmatrix}
   0  & T_{m-\sqrt{-1}n}\\
   T_{m+\sqrt{-1}n} & 0 \\
  \end{smallmatrix}}\right], \left[{\begin{smallmatrix}
   1  & 0\\
   0 & -1 \\
  \end{smallmatrix}}\right]\right)
\]
is a $2^+$-summable even spectral triple on the noncommutative torus $\cla_\theta$. To check nontriviality, one  pairs it with the $K_0$-class of the Powers-Rieffel projection $p_{\theta}$. The projection $p_\theta$  has a power series expression in $u$ and $v$, which we denote by $\mathcal{P}(u,v)$. Now, the following operator
\[
p_\theta\,T_{\frac{m+\sqrt{-1}n}{\sqrt{m^2+n^2}}}p_\theta\,:\,p_\theta \ell^2(\bbz^2)\longrightarrow p_\theta\ell^2(\bbz^2)
\]
is Fredholm and its index is nonzero \cite{Had-2004aa}, which proves nontriviality of this spectral triple.
\smallskip

\textbf{Notation:~} The index of the above Fredholm operator is denoted by the symbol $\varrho$ throughout this section.
\smallskip

Recall Propn. \ref{E_1} from Sec. $(5)$, and let
\begin{IEEEeqnarray}{rCl}\label{clh_k}
 \clh_r=\bigoplus_{i+j=n_r\,,\,k \in \bbz}E_1(i,j,k)\,.
\end{IEEEeqnarray}
Take a unit vector in $E_1(\frac{n_r}{2}, \frac{n_r}{2},n_r)$, and denote it by $|\frac{n_r}{2},\frac{n_r}{2},0\rangle$. 
For $i,j \in \frac{1}{2}\bbz$ such that $i+j=n_r$, we define the following
\begin{IEEEeqnarray}{lCl}\label{ONBelement}
|i,j,k\rangle=\begin{cases}
                    (b^*)^{2i-n_r} (D^*)^{(n_r-2i+k)}|\frac{n_r}{2},\frac{n_r}{2},0\rangle & \mbox{ if } i \geq \frac{n_r}{2}, \cr
                    (b)^{n_r-2i}  (D^*)^k|\frac{n_r}{2}, \frac{n_r}{2},0\rangle & \mbox{ if } i <\frac{n_r}{2}. \cr
                    \end{cases}
\end{IEEEeqnarray}
Observe that 
\begin{IEEEeqnarray}{rCll}\label{act of b}
b\,|i,j,k\rangle &=& \begin{cases}
\big|i-\frac{1}{2},j+\frac{1}{2},k-1\big\rangle & \mbox{ if } i \geq \frac{n_r+1}{2}, \cr
\big|i-\frac{1}{2},j+\frac{1}{2},k\big\rangle & \mbox{ if } i \leq \frac{n_r}{2};
\end{cases}
\end{IEEEeqnarray}
\begin{IEEEeqnarray}{rCll}\label{act of b*}
b^*\,|i,j,k\rangle &=& \begin{cases}
\big|i+\frac{1}{2},j-\frac{1}{2},k+1\big\rangle& \mbox{ if } i \geq \frac{n_r}{2}, \cr
\big|i+\frac{1}{2},j-\frac{1}{2},k\big\rangle& \mbox{ if } i \leq \frac{n_r-1}{2};
\end{cases}
\end{IEEEeqnarray}
and
\begin{IEEEeqnarray}{rCll}\label{act of D}
D\,|i,j,k\rangle &=& e^{2\pi\sqrt{-1}(2i-n_r)\theta}\,|i,j,k-1\rangle\,\,,\quad  & D^*\,|i,j,k\rangle=e^{2\pi\sqrt{-1}(n_r-2i)\theta}\,|i,j,k+1\rangle\,.
\end{IEEEeqnarray}
 Moreover, since $b,b^*,D$ and $D^*$ act as unitary operators when restricted to the invariant space $E_1$, it follows that $|i,j,k\rangle$ is a unit vector in $E_1(i,j,n_r+k)\subseteq\clh_r$. Hence, we can write 
\begin{IEEEeqnarray}{rCl}\label{ONB}
|i,j,k\rangle &=& \sum_{m=0}^{\infty}c_m^{(i,j,k)}e_{i,j, i+j+k+m}^{w_{ij}+m}\,.
\end{IEEEeqnarray}
Note that $\sum_{m=0}^{\infty}|c_m^{(i,j,k)}|^2=1$ as $|i,j,k\rangle $ is a unit vector.

\bppsn\label{c_k}
Define $C_{i,j,k}=\sum_{m=1}^{\infty}\big|c_{m}^{(i,j,k)}\big|^2$ for $k \in \bbz$ and $i,j\in\frac{1}{2}\bbz$ with 
$i+j=n_r \in\Omega$ for some $r\in\bbn$. Then, one has the following.
\begin{enumerate}[$(i)$]
\item $C_{i,j,k}\to 0$ as $i\to\pm\infty$, provided $i+j=n_r$ for some fixed $r\in\bbn$.
\item $C_{i,j,k}\to 0$ as $i+j\to\infty$.
\end{enumerate}
\eppsn
\prf \begin{enumerate}[$(i)$]
      \item Using Thm. \ref{the action} and Eqns. (\ref{act of D},\ref{ONB}), we have
\begin{IEEEeqnarray*}{rCl}
& & D|i,j,k\rangle = D\sum_{m=0}^{\infty}c_m^{(i,j,k)}e_{i,j, i+j+k+m}^{w_{ij}+m} \\
&\Rightarrow& e^{2\sqrt{-1}\pi(2i-n_r)\theta}\,|i,j,k-1\rangle= \sum_{m=0}^{\infty}c_m^{(i,j,k)}D\,e_{i,j, i+j+k+m}^{w_{ij}+m}\\
&\Rightarrow & e^{2\sqrt{-1}\pi(2i-n_r)\theta}\sum_{m=0}^{\infty}c_m^{(i,j,k-1)}e_{i,j, i+j+k+m-1}^{w_{ij}+m} 
= \sum_{m=0}^{\infty}c_m^{(i,j,k)} e^{2\sqrt{-1}\pi(i-j)\theta}\, e_{i,j, i+j+k+m-1}^{w_{ij}+m}\\
\end{IEEEeqnarray*}
Comparing coefficients of both sides, we get $|c_m^{(i,j,k)}|=|c_m^{(i,j,k-1)}|$. This shows that for any $k_1, k_2 \in \bbz$, 
$C_{i,j,k_1}=C_{i,j,k_2}$. Hence, without loss of generality, we can take $k=0$. Fix $r \in \bbn$. 
To get $\lim_{i \rightarrow -\infty}C_{i,j,0}$, 
we assume that $i\leq n_r/2$, so that $w_{ij}=j$.
\begin{IEEEeqnarray*}{rCl}
&&b|i,j,0\rangle = b\sum_{m=0}^{\infty}c_m^{(i,j,0)}e_{i,j, i+j+m}^{j+m} \\
&\Rightarrow& |i-\frac{1}{2},j+\frac{1}{2},0\rangle = 
\sum_{m=0}^{\infty}c_m^{(i,j,0)}\Big(\beta_+^+(j+m,i,j)e_{i+\frac{1}{2}\,,\,j-\frac{1}{2}\,,\,i+j+m+1}^{j+m+\frac{1}{2}}+  \\
&& \hspace{3in}  \beta_-^+(j+m,i,j)e_{i+\frac{1}{2}\,,\,j-\frac{1}{2}\,,\,i+j+m}^{j+m-\frac{1}{2}}\Big)\\
&\Rightarrow & \sum_{m=0}^{\infty}c_m^{(i-\frac{1}{2},j+\frac{1}{2},0)}e_{i-\frac{1}{2},j+\frac{1}{2}, i+j+m}^{j+m+\frac{1}{2}} 
= \sum_{m=0}^{\infty}\Big(c_m^{(i,j,0)}\beta_+^+(j+m,i,j) \\
&& \hspace{2.5in}  +c_{m+1}^{(i,j,0)} 
\beta_+^-(j+m+1,i,j)\Big)
 e_{i+\frac{1}{2}\,,\,j-\frac{1}{2}\,,\,i+j+m+1}^{j+m+\frac{1}{2}}.
\end{IEEEeqnarray*}
Comparing coefficients of both sides, we get
\begin{IEEEeqnarray}{rCl}\label{eqn5}
 c_m^{(i-\frac{1}{2},j+\frac{1}{2},0)}=c_m^{(i,j,0)}\beta_+^+(j+m,i,j)+c_{m+1}^{(i,j,0)} 
\beta_+^-(j+m+1,i,j).
\end{IEEEeqnarray}
By Thm. \ref{the action}, we have the following estimates.
\begin{IEEEeqnarray}{rCl}\label{eqn6}
 |\beta_+^+(\ell,i,j)|\leq |q|^{r-i}, \qquad |\beta_+^-(\ell,i,j)|\leq |q|^{r-j}.
\end{IEEEeqnarray}
From Eqns. (\ref{eqn5},\ref{eqn6}), we get 
\[
 |c_m^{(i-\frac{1}{2},j+\frac{1}{2},0)}|\leq |q|^m|c_m^{(i,j,0)}|  
 +|q|^{m+1}|c_{m+1}^{(i,j,0)}|\leq \frac{|q|^m}{1-|q|}\big(|c_m^{(i,j,0)}|+|q||c_{m+1}^{(i,j,0)}|\big).
\]
By iterating the process starting with $|\frac{n_r}{2},\frac{n_r}{2},0\rangle$, using Eqns. (\ref{eqn5},\ref{eqn6})  
and the fact that $|c_m^{(i,j,m)}|\leq 1$, we have
\[
 |c_m^{(i,j,0)}|\leq \frac{|q|^{m(n_r-2i)}}{1-|q|}\left(|c_m^{(\frac{n_r}{2},\frac{n_r}{2},0)}|+|q| 
 |c_{m+1}^{(\frac{n_r}{2},\frac{n_r}{2},0)}|+ \cdots +|q|^{n_r-2i} |c_{m+n_r-2i}^{(\frac{n_r}{2},\frac{n_r}{2},0)}|\right)\leq 
 \frac{|q|^{m(n_r-2i)}}{(1-|q|)^2}.
\]
This proves that $C_{i,j,0}\rightarrow 0$ as $ i \rightarrow -\infty$. By replacing $b$ with $b^*$ in this computations, one can show that 
$C_{i,j,0}\rightarrow 0$ as $ i \rightarrow \infty$.
\item  
Employing part (i) of the proposition, it is enough to show that $C_{\frac{n_r}{2},\frac{n_r}{2},0} \rightarrow 0$ as $r \rightarrow \infty$.
Take $v \in E_{|q|^{2n_r}}(0,0,0)$ such that $\|v\|=1$. By  part (vii) of Lemma \ref{few observations} and using $\|a\|\leq 1$, we have
\begin{IEEEeqnarray}{rCl}\label{eqn7}
 (a^*)^{n_r}v = M|\frac{n_r}{2},\frac{n_r}{2},0\rangle
\end{IEEEeqnarray}
for some constant $M$ with $|M|\leq 1$. Now, proceeding as above using Eqn. \ref{eqn7} and  the estimates of
$|\alpha_+^+(\ell,i,j)|$ and $\alpha_+^+(\ell,i,j)$ obtained from Thm. \ref{the action} we get  the claim.\qed
\end{enumerate}

\bppsn\label{C_{ij}}
One has the following.
\begin{enumerate}[$(i)$]
\item The set $\{|i,j,k\rangle:k\in\bbz;\,i,j\in\frac{1}{2}\bbz,\,i+j=n_r\}$ is an orthonormal basis of $\clh_r$.
\item The set $\{|i,j,k\rangle:k\in\bbz;\,i,j\in\frac{1}{2}\bbz,\,i+j=n_r,\,r\in\bbn\}$ is an orthonormal basis of $E_1$. 
\end{enumerate}
\eppsn
\prf Observe that $A(i,j,k)$ is orthogonal to $A(i^{'},j^{'},k^{'})$ if $(i,j,k)\neq (i^{'},j^{'},k^{'})$. This along with the decomposition in Eqn. \ref{clh_k} proves the first part. By Propn. \ref{E_1} part $(iv)$, the second part follows.\qed
  
Now, let $P=\mathds{1}_{\{bb^*=1\}}\in C(U_q(2))$. Under the GNS representation $\pi_h,\,P$ is the orthogonal projection onto the subspace $E_1$. Note that $PbP=Pb=bP$ and $PDP=PD=DP$, as both $b$ and $D$ commute with $bb^*$. Consider the $C^*$-algebra $\mathcal{B}$ generated by $Pb$ and $PD$. Then, $\mathcal{B}$ is a unital $C^*$-algebra with the multiplicative unit $P$. Using the defining relations in \ref{relations}, it follows that $PbPD=e^{2\sqrt{-1}\pi\theta}PDPb$, and $PD$ is a unitary in $\mathcal{B}$. Now, the relation $aa^*+bb^*=1$ implies that $Paa^*P+Pbb^*P=P$. It is easy to check that for any vector $v$ in the fixed point subspace under the action of $bb^*$, one has $a^*v=0$ as $\sigma(bb^*)=\{|q|^{2m}:m\in\bbn\}\cup\{0\}$. Thus, it follows that $Pbb^*P=P$ and consequently, $Pb$ is a unitary in $\mathcal{B}$. Hence, by the universality of the noncommutative torus $\cla_{\theta}$ and its simpleness, there exists an isomorphism $\psi: \cla_{\theta}\longrightarrow\mathcal{B}$ sending $u\longmapsto Pb$ and $v\longmapsto PD$. So, one can conclude  that $K_0(\mathcal{B})$ is generated by $[P]$ and $[\mathcal{P}(Pb,PD)]$. Moreover, Lemma (\ref{kernel is NC torus},\ref{needed K groups}) and Thm. \ref{K groups} says that one can view $[p \otimes p_\theta]$, one of the generators of $K_0(C(U_q(2)))$, as $[\mathcal{P}(Pb,PD)]$. We denote $\mathcal{P}(Pb,PD)$ by $P_\theta$. Since $\mathcal{P}(Pb,PD) \in \mathcal{B}$, we get that
\begin{IEEEeqnarray}{lCl}\label{eqn1}
 PP_\theta=P_\theta P=P_\theta\,. 
\end{IEEEeqnarray}
Let $\rho:\mathcal{B}\longrightarrow\mathcal{B}(PL^2(h))$ be a representation of $\mathcal{B}$ given by $\rho(x)=\pi_h(x)$ for $x\in\mathcal{B}$. By Eqns. (\ref{act of b},\ref{act of b*},\ref{act of D}) and Propn. \ref{C_{ij}}, the closed subspace $PL^2(h)$ is invariant under $b,b^*,D$ and $D^*$, and hence the representation $\rho$ is well-defined.
Define an  operator $F:PL^2(h)\longrightarrow PL^2(h)$ by $F\,|i,j,k\rangle=f_0(i,j,k)|i,j,k\rangle$, where 
\begin{IEEEeqnarray}{lCl}\label{F}
f_0(i,j,k)= \begin{cases}
\frac{(2w_{ij}+1)+\sqrt{-1}(k+j-w_{ij})}{\sqrt{(2w_{ij}+1)^2+(k+j-w_{ij})^2}} & \mbox{ if }\,\, i\neq-w_{ij}, \cr
\frac{-(2w_{ij}+1)+\sqrt{-1}(k+j-w_{ij})}{\sqrt{(2w_{ij}+1)^2+(k+j-w_{ij})^2}}& \mbox{ if }\,\, i=-w_{ij}, \cr
\end{cases}
\end{IEEEeqnarray}

\blmma\label{cpt part}
The truncated operator $PT|T|^{-1}P:PL^2(h) \longrightarrow PL^2(h)$ is a compact perturbation of $F$.
\elmma
\prf Let $f(\ell,i,m)=\frac{d(\ell,i,m)}{|d(\ell,i,m)|}$. Then, one has
\begin{IEEEeqnarray*}{rCl}
 T|T|^{-1}\,|i,j,k\rangle &=&T|T|^{-1}\Big(\sum_{m=0}^{\infty}c_m^{(i,j,k)}e_{i,j, i+j+k+m}^{w_{ij}+m}\Big),\\
 &=&\sum_{m=0}^{\infty}c_m^{(i,j,k)}f(w_{ij}+m,i,m+k)e_{i,j, i+j+m+k}^{w_{ij}+m}\,.
\end{IEEEeqnarray*}
Hence, we get that
\begin{IEEEeqnarray*}{rCl}
  PT|T|^{-1}P\,|i,j,k\rangle &=& \sum \Big\langle T|T|^{-1}\,|i,j,k\rangle,\,|i^{'},j^{'},k^{'}\rangle \Big\rangle\,|i^{'},j^{'},k^{'}\rangle\\
  &=&\Big\langle T|T|^{-1}\,|i,j,k\rangle,\,|i,j,k\rangle\Big\rangle\,|i,j,k\rangle \\
  &=& \Big(\sum_{m=0}^{\infty} \big{|}c_m^{(i,j,k)}\big{|}^2f(w_{ij}+m,i,i+j+k+m)\Big)\,|i,j,k\rangle\,.
\end{IEEEeqnarray*}
Using the facts that $\sum_{m=0}^{\infty} \big{|}c_m^{(i,j,k)}\big{|}^2=1$ and $f(w_{ij}+m,i,i+j+k+m)=f_0(i,j,k)$ one has the following,
\begin{IEEEeqnarray*}{rCl}
(PT|T|^{-1}P-F)\,|i,j,k\rangle
&=& \Big(\sum_{m=1}^{\infty} \big{|}c_m^{(i,j,k)}\big{|}^2\big(f(w_{ij}+m,i,i+j+k+m)-f_0(i,j,k)\big)\Big)\,|i,j,k\rangle\,.\\
\end{IEEEeqnarray*}
Now, the claim follows from Propn. \ref{c_k} along with the following observations,
\[
|f(w_{ij}+m,i,i+j+k+m)-f_0(i,j,k)|\leq 2\,,
\]
and
\[
|f(w_{ij}+m,i,i+j+k+m)-f_0(i,j,k)| \to 0\,\, \mbox{ as } m\to 0\,.
\]
\qed

We now decompose the operator $P_{\theta}FP_{\theta}$ in the following manner. 
For $r\geq 0$, let $P_l$ be the projection onto the subspace $\clh_r$ of $PL^2(h)$. 
It follows from Eqns. (\ref{act of b},\ref{act of D}) that $\clh_r$ is invariant under the action of $\mathcal{B}$, i.e, $P_rx=xP_r$ for all $x \in\mathcal{B}$.   
This induces a representation $\rho_r: \mathcal{B}\longrightarrow \mathcal{B}(\clh_r)$ 
defined by $\rho_r(x)=\rho(x)\restr{\clh_r}$ for $x \in\mathcal{B}$. For $r \geq 0$, let $F_r= F\restr{\clh_r}$  
and $P_{\theta}^r=(P_{\theta})\restr{\clh_r}$. Hence, we have the following,
\begin{IEEEeqnarray}{lCl}\label{directsum}
\rho=\oplus_{r=0}^{\infty}\,\rho_r\,\,,\quad P=\oplus_{r=0}^{\infty}\,P_r\,\,,\quad F=\oplus_{r=0}^{\infty}\,F_r\,\,,\quad P_{\theta}FP_{\theta}=\oplus_{r=0}^{\infty}\,P_{\theta}^rF_rP_{\theta}^r\,.
\end{IEEEeqnarray}

\blmma\label{ind_0}
The tuple 
$\big(\mathcal{B}, \rho_0 \oplus \rho_0,\left[ {\begin{smallmatrix}
   0 & F_0^*\\
   F_0 & 0 \\
  \end{smallmatrix}} \right],\left[ {\begin{smallmatrix}
   1  & 0\\
   0 & -1 \\
  \end{smallmatrix}} \right]\big)$ is a finitely summable even Fredholm module, and $\,\mathrm{ind}(P_{\theta}^0F_0P_{\theta}^0)=\varrho\,$.
\elmma
\prf Define a unitary operator $W_0: \clh_0 \longrightarrow \ell^2(\bbz^2)$ as follows~:
\begin{IEEEeqnarray}{lCl}\label{W_0}
W_0\,|i,j,k\rangle &=& 
e_{2i} \otimes e_k\,\,, \quad \mbox{ for }\,k \in \bbz, i,j \in \frac{1}{2}\bbz, i+j=0.
\end{IEEEeqnarray}
The tuple $\big(\mathcal{B}, \rho_0\oplus \rho_0,\left[ {\begin{smallmatrix}
   0 & F_0^*\\
   F_0 & 0 \\
  \end{smallmatrix}} \right],\left[ {\begin{smallmatrix}
   1  & 0\\
   0 & -1 \\
  \end{smallmatrix}} \right]\big)$ is unitary equivalent to the following finitely summable Fredholm module
\[
\Big(\mathcal{B}, \tau \circ \psi^{-1} \oplus \tau \circ \psi^{-1}, \left[ {\begin{smallmatrix}
   0 & T_{\frac{t-1-\sqrt{-1}s}{(t-1)^2+s^2}}\\
  T_{\frac{t-1+\sqrt{-1}s}{(t-1)^2+s^2}} & 0 \\
  \end{smallmatrix}} \right],\left[ {\begin{smallmatrix}
   1  & 0\\
   0 & -1 \\
  \end{smallmatrix}} \right]\Big)
  \]
through the unitary $W_0$, which is  homotopic to the Fredholm module
\[
\Big(\mathcal{B}, \tau \circ \psi^{-1} \oplus \tau \circ \psi^{-1}, \left[ {\begin{smallmatrix}
   0 & T_{\frac{t-\sqrt{-1}s}{t^2+s^2}}\\
  T_{\frac{t+\sqrt{-1}s}{t^2+s^2}} & 0 \\
  \end{smallmatrix}} \right],\left[ {\begin{smallmatrix}
   1  & 0\\
   0 & -1 \\
  \end{smallmatrix}} \right]\Big)\,.
  \]
This completes the proof.\qed

\blmma \label{ind_k}
For $r > 0$, the tuple 
$\zeta_r=\big(\mathcal{B}, \rho_r \oplus \rho_r,\left[ {\begin{smallmatrix}
   0 & F_r^*\\
   F_r & 0 \\
  \end{smallmatrix}} \right],\left[ {\begin{smallmatrix}
   1  & 0\\
   0 & -1 \\
  \end{smallmatrix}} \right]\big)$ is  a  finitely summable  even Fredholm module. Moreover, for $r_1,r_2>0$, one has
\[
 ind(P_{\theta}^{r_1}F_{r_1}P_{\theta}^{r_1})=ind(P_{\theta}^{r_2}F_{r_2}P_{\theta}^{r_2}).
\]
\elmma
\prf For $i,j \in \frac{1}{2}\bbz, i+j=n_r$, one has
\begin{IEEEeqnarray*}{rCl}
  [F_r,Pb]\,|i,j,k\rangle &=& F_rPb\,|i,j,k\rangle-PbF\,|i,j,k\rangle\\
  &=& F_r\,|i-1/2,j+1/2, k\rangle-f(i,j,k)Pb\,|i,j,k\rangle\\
  &=& \big(f(i-1/2,j+1/2, k)-f(i,j,k)\big)\,|i-1/2,j+1/2, k\rangle
\end{IEEEeqnarray*}
Using eqn. (\ref{F}), it is straightforward to check that $[F_r,Pb]$ is a compact operator on $\clh_r$. Similarly one can show that $[F_r,PD]$ is compact, which proves that $\zeta_r$ is a finitely summable even Fredholm module.

Fix $r >0$. To prove the last part, we will show that the Fredholm module $\zeta_r$ is unitary equivalent to a Fredholm module which is homotopic to $\zeta_1$. Define a unitary operator $W_r: \clh_r \longrightarrow \clh_1$ by the following,
\begin{IEEEeqnarray}{lCl}
W_r\,|i,j,k\rangle &=& 
|i+n_1/2-n_r/2,j+n_1/2-n_r/2,k\rangle\,,\,\,\mbox{for } k \in \bbz, i,j \in \frac{1}{2}\bbz, i+j=n_r.
\end{IEEEeqnarray}
It follows from Eqns. (\ref{act of b},\ref{act of D}) that $W_r\rho_r(b)W_r^*=\rho_1(b)$ and $W_r\rho_r(D)W_r^*=\rho_1(D)$, which further implies that $W_r\rho_rW_r^*=\rho_1$. Therefore, we have
\begin{IEEEeqnarray}{lCl}\label{P}
 W_rP_{\theta}^rW_r^*=P_{\theta}^{1}\,.
\end{IEEEeqnarray}
For $0\leq t \leq 1$, define
\[
 F_r^t\,|i,j,k\rangle=\frac{2w_{ij}+1+t(n_r-n_1)+\sqrt{-1}(k+j-w_{ij})}{\sqrt{(2w_{ij}+1+t(n_r-n_1))^2+(k+j-w_{ij})^2}}\,|i,j,k\rangle\,. 
\]
Consider  the tuple  $\,\zeta_r^t=\left(B, \rho_1\oplus \rho_1,\left[ {\begin{smallmatrix}
   0 & (F_r^t)^*\\
   F_r^t & 0 \\
  \end{smallmatrix}} \right],\left[ {\begin{smallmatrix}
   1  & 0\\
   0 & -1 \\
  \end{smallmatrix}} \right]\right)$.
Similar to the case of $\zeta_r$, one can show that $\zeta_r^t$ is a finitely summable even Fredholm module.  Moreover, $\zeta_0^r=\zeta_1$ and $\zeta_1^r=W_r\zeta_rW_r^*$. Thus, $\zeta_1$ is homotopic to   $W_r\zeta_rW_r^*$ and hence, they represent the same element in $K^0$-group. Using the $K_0-K^0$ pairing, one gets that
\begin{IEEEeqnarray*}{lCl}
 ind(P_{\theta}^{1}F_{1}P_{\theta}^{1}) &=& ind(P_{\theta}^{1}W_rF_1W_r^*P_{\theta}^{1})\\
 &=& ind(W_rP_{\theta}^rF_{1}P_{\theta}^rW_r^*) \qquad (\mbox{ by eqn. }\ref{P})\\
 &=& ind(P_{\theta}^rF_1P_{\theta}^r)\,,
\end{IEEEeqnarray*}
which completes the proof.\qed

\bppsn\label{index_0}
The operator $P_{\theta}FP_{\theta}:P_{\theta}L^2(h)\longrightarrow P_{\theta}L^2(h)$ is a Fredholm operator 
with $ind(P_{\theta}FP_{\theta})=\varrho$.
\eppsn
\prf The tuple 
$(\mathcal{B}, \rho \oplus \rho,\left[ {\begin{smallmatrix}
   0 & F^*\\
   F & 0 \\
  \end{smallmatrix}} \right],\left[ {\begin{smallmatrix}
   1  & 0\\
   0 & -1 \\
  \end{smallmatrix}} \right])$
is a finitely summable even Fredholm module as 
\[
\big(\mathcal{B}, \rho \oplus \rho,\left[ {\begin{smallmatrix}
   0 & F^*\\
   F & 0 \\
  \end{smallmatrix}} \right],\left[ {\begin{smallmatrix}
   1  & 0\\
   0 & -1 \\
  \end{smallmatrix}}\right]\big)=\bigoplus_{r=0}^{\infty}\,\big(\mathcal{B}, \rho_r \oplus \rho_r,\left[{\begin{smallmatrix}
   0 & F_r^*\\
   F_r & 0 \\
  \end{smallmatrix}}\right],\left[{\begin{smallmatrix}
   1  & 0\\
   0 & -1 \\
  \end{smallmatrix}}\right]\big)\]
This shows that the operator $P_{\theta}FP_{\theta}:P_{\theta}L^2(h) \longrightarrow P_{\theta}L^2(h)$ is a Fredholm operator.  Therefore, we have $ind(P_{\theta}FP_{\theta})<\infty$. By Eqn. \ref{directsum}, one has 
\[
 ind(P_{\theta}FP_{\theta})=\sum_{r=0}^{\infty}ind(P_{\theta}^rF_rP_{\theta}^r). 
\]
It follows from Lemma \ref{ind_k} that $ind(P_{\theta}^rF_rP_{\theta}^r)=0$ for $r\geq 1$. Now, by Lemma \ref{ind_0}, one can conclude that $\mathrm{ind}(P_{\theta}FP_{\theta})=\varrho$.\qed
\smallskip

\textbf{Proof of Thm. \ref{final thm}~:} Observe that $P_{\theta}T|T|^{-1}P_{\theta}=P_{\theta}PT|T|^{-1}PP_{\theta}$ as $P_{\theta}P=PP_{\theta}=P_{\theta}$. Therefore by Lemma \ref{cpt part}, the operator $P_{\theta}T|T|^{-1}P_{\theta}$ is a compact perturbation of $P_{\theta}FP_{\theta}$. From Propn. \ref{index_0}, we get that
\[
ind(P_{\theta}T|T|^{-1}P_{\theta})=ind(P_{\theta}FP_{\theta})=\varrho. 
\]
Thus, the value of the $K_0-K^0$ pairing coming from the Kasparov product is given by the following,
\[
 \langle [p \otimes p_{\theta}], (\mathscr{A}_q\,,\pi_{eq}\,,\mathscr{D}\,,\gamma)\rangle =ind(P_{\theta}T|T|^{-1}P_{\theta})=\varrho\,,
\]
and this completes the proof.\qed
\smallskip

We finally conclude the case of $\theta=0,1$. In these two situations, $q$ is real. Recall from Sec. \ref{Sec 2} that in these cases $K_0(C(U_q(2)))$ is generated by $[1]$ and $[p\otimes Bott]$, where ``$Bott$'' denotes the Bott projection is $M_2(C(\mathbb{T}^2))$. It is well known that the following spectral triple
\[
\left(C^\infty(\mathbb{T}^2)\,,\,\ell^2(\bbz^2)\otimes\mathbb{C}^2,\left[ {\begin{smallmatrix}
   0  & \frac{\partial}{\partial x}-\sqrt{-1}\frac{\partial}{\partial y}\\
   \frac{\partial}{\partial x}+\sqrt{-1}\frac{\partial}{\partial y} & 0 \\
  \end{smallmatrix} } \right], \left[ {\begin{smallmatrix}
   1  & 0\\
   0 & -1 \\
  \end{smallmatrix} } \right]\right)
\]
is nontrivial. The $C^*$-algebra $\mathcal{B}$ becomes $C(\mathbb{T}^2)$ and all the proofs in this section remain valid, and hence we have that Thm. \ref{final thm} holds for $\theta=0,1$ also.

%%%%%%%%%%%%%%%%%%%%%%%%%%%%%%%%%%%%%%%
%%%%%%%%%%   Spectral dimension   %%%%%%%%%%%%%%%%%%%
%%%%%%%%%%%%%%%%%%%%%%%%%%%%%%%%%%%%%%%

\newsection{The spectral dimension}\label{Sec 8}

Spectral dimension of a compact quantum group has been introduced in \cite{CP-2016aa}. In this section we compute it for the compact quantum group $U_q(2)$. We will mainly follow \cite{Sau-2020aa} for the computation. To put into the appropriate framework we define the following,
\begin{IEEEeqnarray*}{rCl}
 \Gamma &=& \{\gamma=(\gamma_1,\gamma_2,\gamma_3): \gamma_1 \in \frac{1}{2}\bbn, \gamma_2 \in \bbz,
 \gamma_3 \in \{-\gamma_1,- \gamma_1+1, \cdots  \gamma_1\}\} \\
 e^{\gamma} &=& a^{\gamma_1-\gamma_3}b^{\gamma_1+\gamma_3}D^{\gamma_2}=t_{\gamma_3,-\gamma_1}^{\gamma_1}D^{\gamma_2},\\
 e^{(\gamma,j)} &=& t_{\gamma_3,j}^{\gamma_1}D^{\gamma_2}\,;\,-\gamma_1 \leq j \leq \gamma_1,\\
 \epsilon_1&=&(1/2,0,0)\,, \qquad \epsilon_2=(0,1,0)\,, \qquad  \epsilon_3=(0,0,1/2). 
\end{IEEEeqnarray*}

\blmma \label{bound hwv}
With the above set up, one has the followings.
\begin{enumerate}[$(i)$]
\item
 $\sup_{\{\gamma \in \Gamma\}}\frac{\|e^{\gamma}\|}{\|e^{\gamma+\epsilon_2}\|} <\infty$.
 \item
 $\sup_{\{\gamma \in \Gamma\}}\frac{\|e^{\gamma}\|}{\|e^{\gamma+\epsilon_1-\epsilon_3}\|} <\infty$. 
 \item
 $\sup_{\{\gamma \in \Gamma\}}\frac{\|e^{\gamma}\|}{\|e^{\gamma+\epsilon_1+\epsilon_3}\|} <\infty$. 
\end{enumerate}
\elmma 
\prf  Note  that for $\gamma=(\gamma_1,\gamma_2,\gamma_3) \in \Gamma$, 
$e^{\gamma}=t_{\gamma_3,-\gamma_1}^{\gamma_1}D^{\gamma_2}$. Hence by Thm $4.17$ in \cite{GuiSau-2020aa}, 
we have $\|e^{\gamma}\|=\frac{|q|^{\gamma_3}}{\sqrt{|2\gamma_1+1|_{|q|}}}$. Applying this result, we get 
 $\frac{\|e^{\gamma}\|}{\|e^{\gamma+\epsilon_2}\|}=1$. Moreover, 
\begin{IEEEeqnarray*}{lCl}
 \frac{\|e^{\gamma}\|}{\|e^{\gamma+\epsilon_1-\epsilon_3}\|}&=& 
 \frac{|q|^{2\gamma_3}}{|q|^{2\gamma_3-1}}\sqrt{\frac{|2\gamma_1+2|_{|q|}}{|2\gamma_1+1|_{|q|}}}
 = |q|\sqrt{\frac{|q|^{2\gamma_1}(1+|q|^2+\cdots + |q|^{4\gamma_1+2})}{|q|^{2\gamma_1+1}(1+|q|^2+\cdots + |q|^{4\gamma_1})}}\\
 &\leq &  \sqrt{2|q|}.
\end{IEEEeqnarray*}
Similarly, 
\begin{IEEEeqnarray*}{lCl}
 \frac{\|e^{\gamma}\|}{\|e^{\gamma+\epsilon_1+\epsilon_3}\|}&=& 
 \frac{|q|^{2\gamma_3}}{|q|^{2\gamma_3+1}}\sqrt{\frac{|2\gamma_1+2|_{|q|}}{|2\gamma_1+1|_{|q|}}}
 = \frac{1}{|q|}\sqrt{\frac{|q|^{2\gamma_1}(1+|q|^2+\cdots + |q|^{4\gamma_1+2})}{|q|^{2\gamma_1+1}(1+|q|^2+\cdots + |q|^{4\gamma_1})}}\\
 &\leq &  \frac{\sqrt{2}}{|q|}\,.
\end{IEEEeqnarray*}
\qed

Let $c$ be an upper bound of all the suprema mentioned in Lemma \ref{bound hwv}. Take $R=\{a,b,D,D^*\}$.
From \cite{Sau-2020aa}, recall the growth graph $\mathcal{G}_R^c\,$. One takes the vertex set  of $\mathcal{G}_R^c$ to be $\Gamma$. For $x \in R$, we write  $\gamma \rightsquigarrow_{x} \gamma^{'}$ if 
\begin{IEEEeqnarray}{rCl}\label{eq1}
e^{(\gamma^{'},j^{'})}=re^{(\gamma,j)}\,\,\,\mbox{ and } \,\,\,\frac{\|  e^{(\gamma,j)}\|}{\|e^{(\gamma^{'},j^{'})}\|}<c
\end{IEEEeqnarray}
for some $ 1\leq j \leq N_{\gamma}$ and $ 1\leq j^{'} \leq N_{\gamma^{'}}$. Define the edge set of $\mathcal{G}_R^c$ to be the following,
\[
 E:=\big\{(\gamma,\gamma^{'}): \gamma \rightsquigarrow_{y} \gamma^{'}   \mbox{ for some } y \in R \big\}.
\]
We write $\gamma\to\gamma^{'}$ if $(\gamma,\gamma^{'})\in E$. We say that the graph $\mathcal{G}_R^c$ has a root $\gamma_0 \in \Gamma$ if for any $\gamma \in \Gamma$, there is a directed path from $\gamma_0$ to $\gamma$. The following Lemma says that $\mathcal{G}_R^c$ has a root $(0,0,0)$.

\blmma\label{path}
For $\gamma=(\gamma_1,\gamma_2,\gamma_3) \in \Gamma$,  there is a directed path in $\mathcal{G}_R^c$  from $(0,0,0)$ to $\gamma$, and it is of length less than or equal to $2\gamma_1+|\gamma_2|$.
\elmma
\prf  By Lemma \ref{bound hwv}, we have 
\[
 \gamma \rightsquigarrow_{D} \gamma + \epsilon_2\,\,,\qquad  \gamma \rightsquigarrow_{D^*} \gamma - \epsilon_2\,\,,\qquad \gamma \rightsquigarrow_{a} \gamma + \epsilon_1- \epsilon_3\,,
\]
for $\gamma \in \Gamma$. For $\gamma\in \Gamma$ such that $\gamma_1=\gamma_3$, we have $\gamma \rightsquigarrow_{b} \gamma + \epsilon_1+\epsilon_3$. 
Take $\gamma=(\gamma_1,\gamma_2,\gamma_3) \in \Gamma$. One has
\[
 (0,0,0) \rightsquigarrow_{b} (1/2,0,1/2) \rightsquigarrow_{b} \cdots\cdots \rightsquigarrow_{b} 
 \left(\frac{\gamma_1+\gamma_3}{2},0, \frac{\gamma_1+\gamma_3}{2}\right).
\]
Moreover, we have 
\[
 \left(\frac{\gamma_1+\gamma_3}{2},0, \frac{\gamma_1+\gamma_3}{2}\right)\rightsquigarrow_{a}
 \left(\frac{\gamma_1+\gamma_3}{2}+1/2,0, \frac{\gamma_+\gamma_3}{2}-1/2\right) 
 \rightsquigarrow_{a} \cdots\cdots \rightsquigarrow_{a} (\gamma_1,0,\gamma_3).
\]
If $\gamma_2 \geq 0$, we have 
\[
 (\gamma_1,0,\gamma_3)\rightsquigarrow_{D}  (\gamma_1,1,\gamma_3)\rightsquigarrow_{D} \cdots\cdots  \rightsquigarrow_{D} (\gamma_1,\gamma_2,\gamma_3).
\]
If $\gamma_2\leq 0$, we have 
\[
 (\gamma_1,0,\gamma_3)\rightsquigarrow_{D^*}  (\gamma_1,-1,\gamma_3)\rightsquigarrow_{D^*} \cdots\cdots  \rightsquigarrow_{D^*} 
 (\gamma_1,\gamma_2,\gamma_3).
\]
This proves the claim, as the length of this path is $2\gamma_1+|\gamma_2|$.\qed

\bppsn\label{trivialDirac}
Let $L:L^2(h)\longrightarrow L^2(h)$ be an unbounded self-adjoint operator with dense domain $\mathscr{A}_q$ given by the following,
\[
 L\left(e^{(\gamma,j)}\right)= (2\gamma_1+|\gamma_2|)e^{(\gamma,j)}, \mbox{ for } 
 \gamma=(\gamma_1,\gamma_2,\gamma_3) \in \Gamma \mbox{  and } -\gamma_1 \leq j \leq \gamma_1. 
\]
Then, the tuple $\big(\mathscr{A}_q,L^2(h),L\big)$ is a $4^+$-summable spectral triple equivariant under the comultiplication action of $U_q(2)$.
\eppsn
\prf Clearly, the operator $L$ is an unbounded self-adjoint operator with compact resolvent equivariant under the $U_q(2)$ action. Take any $x\in\{a,b,D\}$. By 
 Thm. in \cite{GuiSau-2020aa}, one can see that 
 \[
  xe^{(\gamma,j)}= \sum_{\substack{\beta \in \Gamma \\ \gamma_1-1/2 \leq \beta_1 \leq \gamma_1+1/2 \\ 
  \gamma_{2}-1 \leq \beta_2 \leq \gamma_{2}+1 \\ \gamma_3-1/2 \leq \beta_3 \leq \gamma_3+1/2}}c_{(\beta,j)}e^{(\beta,j)}.
\]
Using this, we get
\begin{IEEEeqnarray*}{rCl}
 \|[L, x]e^{(\gamma,j)}\|
 &\leq& \|xe^{(\gamma,j)}\| \leq \|x\|\|e^{(\gamma,j)}\|\leq \|e^{(\gamma,j)}\|.
\end{IEEEeqnarray*}
Appropriate summability  of the  spectral triple $(\mathscr{A}_q,L^2(h),L)$  immediately follows from the definition of $L$. \qed
\bthm The spectral dimension of $U_q(2)$ is $4$.
\ethm
\prf let $\gamma_0=(0,0,0) \in \Gamma$. Define the length function $\ell_{\gamma_0}:\Gamma \longrightarrow \bbn $ as follows;
\begin{IEEEeqnarray*}{rCl}
 \ell_{\gamma_0}(\gamma)= \begin{cases}
                           1 & ;\,\mbox{ if } \gamma =\gamma_0,\cr
                           \mbox{ length of a shortest path from } \gamma_0 \mbox{ to } \gamma & ;\,\mbox{ otherwise.} \cr
                          \end{cases}
\end{IEEEeqnarray*}
Let $L_{\gamma_0}$ be the unbounded positive operator acting on $L^2(h)$ with dense domain 
$\mathscr{A}_q$ sending  $e^{(\gamma,j)}$ to $\ell_{\gamma_0}(\gamma)e^{(\gamma,j)}$ for all $-\gamma_1 \leq j \leq \gamma_1$ 
and $\gamma \in \Gamma$. It follows from Lemma \ref{path} that $\ell_{\gamma_0}(\gamma)\leq 2\gamma_1+|\gamma_2|$. Therefore, we get that
\[
 \inf\big{\{}p:\,\mbox{Tr}(L_{\gamma_0}^{-p})<\infty\big{\}}\geq 4.
\]
By Propn. $2.4$ in \cite{Sau-2020aa}, it follows that spectral dimension of $U_q(2)$ is greater than or equal to $4$. This result  together with Propn. \ref{trivialDirac} proves the claim.\qed
\bigskip

\section*{Acknowledgements}
Satyajit Guin acknowledges the support of SERB grant MTR/2021/000818 and DST INSPIRE Faculty grant DST/INSPIRE/04/2015/000901. Bipul Saurabh acknowledges the support of SERB grant SRG/2020/000252 and NBHM grant $02011/19/2021/\mbox{NBHM(R.P)/R\&DII}/8758$.
\bigskip

%%%%%%%%%%%%%%%%%%%%%%%%%%%%%%%%%%%%%%%
%%%%%%%%%   References   %%%%%%%%%%%%%%%%%%%%%%%
%%%%%%%%%%%%%%%%%%%%%%%%%%%%%%%%%%%%%%%

\bigskip

\bigskip

\noindent{\sc Satyajit Guin} (\texttt{sguin@iitk.ac.in})\\
         {\footnotesize Department of Mathematics and Statistics,\\
         Indian Institute of Technology, Kanpur,\\
         Uttar Pradesh 208016, India}
\bigskip

\noindent{\sc Bipul Saurabh} (\texttt{bipul.saurabh@iitgn.ac.in})\\
         {\footnotesize Indian Institute of Technology, Gandhinagar,\\  Palaj, Gandhinagar 382355, India}

%%%%%%%%%%%%%%%%%%%%%%%%%%%%%%%%%%%%%%%%%%%%%%%%%%%%%%

\end{document}